\documentclass[10pt]{article}

\title{Rational Homotopy Theory and Differential Graded Category}

\date{Department of Mathematics, Faculty of Science, 
Kyoto University.}

\author{Syunji Moriya\footnote{Corresponding address: Department of Mathematics, Faculty of Science, 
Kyoto University, 
Kyoto, 606-8502, Japan.          \   
E-mail adress: \texttt{moriyasy@math.kyoto-u.ac.jp} \ 
Telephone number: 81-075-753-3700 \ 
FAX number: 81-075-753-3711}}

\usepackage{amsmath,amsthm,amssymb,amscd,mathrsfs}
\usepackage[arrow,matrix]{xy}

\xyoption{all}
\usepackage{enumerate}

\theoremstyle{plain}

\newtheorem{defi}{Definition}[subsection]

\newtheorem{prop}[defi]{Proposition}

\newtheorem{lem}[defi]{Lemma}

\newtheorem{thm}[defi]{Theorem}

\newtheorem{cor}[defi]{Corollary}

\newtheorem{exa}[defi]{Example}

\newcommand{\mf}[1]{{\mathfrak{#1}}}
\newcommand{\mb}[1]{{\mathbf{#1}}}
\newcommand{\bb}[1]{{\mathbb{#1}}}
\newcommand{\mca}[1]{{\mathcal{#1}}}
\newcommand{\ms}[1]{{\mathsf{#1}}}
\newcommand{\msc}[1]{\mathscr{#1}}

\newcommand{\sset}{\ms{sSet}}
\newcommand{\ssetp}{\ms{sSet}_*}
\newcommand{\ssetpc}{\ms{sSet}_*^c}
\newcommand{\ssetpf}{\ms{sSet}_*^f}
\newcommand{\ssetpfq}{\ms{sSet}_*^{f\bb{Q}}}
\newcommand{\dgcu}{\ms{dgCat}}
\newcommand{\dgc}{\ms{dgCat}^{\geq 0}}

\newcommand{\dgcl}{\ms{dgCat}^{\geq 0}_{cl}}

\newcommand{\dgclp}{\ms{dgCat}^{\geq 0}_{cl,*}}

\newcommand{\tannf}{\ms{Tan}^f}
\newcommand{\dggr}{\ms{dgGrph}^{\geq 0}}

\newcommand{\aut}{\mathit{Aut}^{\otimes}}

\newcommand{\modg}{\ms{Mod}(G)}
\newcommand{\dgmodg}{\ms{dgMod}(G)}

\newcommand{\dgalgg}{\ms{dgAlg}(G)}
\newcommand{\dgalgh}{\ms{dgAlg}(H)}
\newcommand{\eqalg}{\ms{EqdgAlg}}
\newcommand{\eqalgf}{\ms{EqdgAlg}^f_1}
\newcommand{\eqalgfp}{\ms{EqdgAlg}^f_{1,*}}
\newcommand{\repg}{\mathit{Rep}(G)}
\newcommand{\reph}{\mathit{Rep}(H)}
\newcommand{\lock}{\mathit{Loc}(K)}
\newcommand{\vect}{\mathit{Vect}}

\newcommand{\ho}{\ms{Ho}}

\newcommand{\fcat}{\mca{F}_{\mathit{cat}}}

\newcommand{\fcl}{\mca{F}_{cl}}
\newcommand{\ob}{\mathit{Ob}}
\newcommand{\mor}{\mathit{Mor}}
\newcommand{\homo}{\mathit{Hom}}
\newcommand{\inhom}{\mf{Hom}}
\newcommand{\inh}{\mf{H}}

\newcommand{\colim}{\mathrm{colim}}

\newcommand{\sob}{S_{\mathit{ob}}}
\newcommand{\smor}{S_{\mathit{mor}}}

\newcommand{\srela}{S_{\mathit{rel}}}
\newcommand{\vtx}{\mathit{Ver}}
\newcommand{\ed}{\mathit{Ed}}
\newcommand{\tens}{\mathit{Tens}}

\newcommand{\ass}{\mathit{Ass}}
\newcommand{\uni}{\mathit{Unit}}
\newcommand{\comm}{\mathit{Comm}}
\newcommand{\Ev}{\mathit{Ev}}
\newcommand{\Coev}{\mathit{Coev}}
\newcommand{\Int}{\mathit{Int}}
\newcommand{\oli}{\overline}

\newcommand{\ev}{\mathit{ev}}
\newcommand{\coev}{\mathit{coev}}
\newcommand{\tdr}{T_{dR}}
\newcommand{\tdrk}{T_{dR}(K)}
\newcommand{\tdrl}{T_{dR}(L)}

\newcommand{\rsp}{\bb{R}\mathit{Sp}}
\newcommand{\loc}{\msc{L}}
\newcommand{\ot}{\! \otimes \!}

\newcommand{\wcl}{W_{cl}}
\newcommand{\T}{\ms{T}}

\newcommand{\bt}{\boxtimes}
\newcommand{\at}{\nabla (*,*)}
\newcommand{\atp}{\nabla (p,*)}
\newcommand{\atq}{\nabla (*,q)}
\newcommand{\vgr}{V^G_r}
\newcommand{\vhr}{V^H_r}

\newcommand{\tF}{\widetilde F}

\begin{document}
\setlength{\baselineskip}{15pt}

\maketitle

\begin{abstract}

We propose a generalization of Sullivan's de Rham homotopy theory to non-simply connected spaces. The formulation is such that the real homotopy type of a manifold should be the closed tensor dg-category of flat bundles on it much the same as the real homotopy type of a simply connected manifold is the de Rham algebra in original Sullivan's theory. We prove the existence of a model category structure on the category of small closed tensor dg-categories and as a most simple case, confirm an equivalence between the homotopy category of spaces whose fundamental groups are finite and whose higher homotopy groups are finite dimensional rational vector spaces and the homotopy category of small closed tensor dg-categories satisfying certain conditions.

\end{abstract}

\begin{center}
\textit{Keywords: rational homotopy theory, non-simply connected space,\\ model category, dg-category.} 
\end{center}
\begin{center}
\end{center}

\section{Introduction}
A rationalization of a simply connected space  $X$ is a map $f:X\to X_{\bb{Q}}$ such that the higher homotopy groups of $X_{\bb{Q}}$ are uniquely divisible and  $f$ induces an isomorphism $\pi _n(X)\otimes _{\bb{Z}}\bb{Q}\cong \pi_n(X_{\bb{Q}})$ for each $n\geq 2$.  We call the homotopy type of $X_{\bb{Q}}$ the rational homotopy type of $X$ and say $X$ is rational if $X_{\bb{Q}}$ is homotopy equivalent to $X$. Sullivan showed  rational homotopy type of a simply connected space of finite type can be recovered from its polynominal de Rham algebra and  the homotopy category of simply connected rational 
spaces of finite $\bb{Q}$-type are equivalent to the homotopy category of 1-connected commutative dg-$\bb{Q}$-algebras of finite type (see 
\cite{sul} or \cite{gug} where the authors call the equivalence the Sullivan-de Rham equivalence).
A feature of Sullivan's theory is that if one consider a $C^{\infty }$-manifold, 
the corresponding dg-algebra over real numbers  is 
(quasi-isomorphic to) the de Rham algebra of the manifold. Because of this feature,  Sullivan's theory has 
geometric applications. See \cite{suldel,sul,mor}.\\
\indent  In the non-simply connected case, as a genaralization of rationalization, Bousfield and Kan \cite{kan} constructed a fiberwise rationalization. For a possibly non-simply connected space $X$, a fiberwise rationalization is a map $f:X\to X_{\bb{Q}}$ such that it induces an isomorphism of fundamental groups and the  map $\tilde f:\widetilde X\to \widetilde {X_{\bb{Q}}}$ between universal coverings is a rationalization in the above sense. We call the homotopy type of $X_{\bb{Q}}$ the rational homotopy type of $X$. For this notion, A.G$\acute{\textrm{o}}$mez-Tato, S.Halperin and D.Tanr$\acute{\textrm{e}}$ \cite{nsrat} generalized the Sullivan's result to non-simply connected spaces. They proposed the notion of 
local system of commutative cochain algebras as  a generalization of commutative dg-algebra  and prove rational homotopy type of spaces with finite rank homotopy groups can be recovered from the corresponding local system and prove an equivalence theorem for non-simply connected 
rational spaces with $\bb{Q}$-finite dimensional higher homotopy groups. As other generalizations of the Sullivan's theory for non-simply connected spaces,  in \cite{trian} a more rigid notion of rational homotopy type is studied and in \cite{brownszc} equivariant dg-algebras are used as algebraic models.  \\

\indent In this paper we introduce  different algebraic  object  viewed as a generalization of commutative dg-algebra
 and as a first step, prove the category of the algebraic objects admits a model category structure (Theorem \ref{result1}) and prove an equivalence theorem for rational spaces 
whose fundamental groups are finite (Theorem \ref{result2}). 
The algebraic objects are small closed tensor differential graded (dg-) categories.\\
\indent A closed tensor dg-category is, roughly speaking, a dg-category which is equipped with a structure of closed symmetric monoidal category consistent with the differential graded structure (see Definition \ref{defofcldgc}).  If one views a dg-algebra as a dg-category with only one object, a symmetric monoidal structure on a dg-category is a natural generalization of commutativity of  dg-algebra and  we need to consider closedness of the symmetric monoidal structure. We mainly consider the pointed case and the corresponding augmented objects are closed tensor dg-categories with fiber functors. Here, a fiber functor of  a closed tensor dg-category $C$ is a dg-functor from $C$ to $\vect$, the closed tensor  category of $\bb{Q}$-vector spaces, which preserves closed tensor structures (see Definition \ref{defofdgclp}).\\
\indent A feature of our formulation is that if one consider over the real (or complex) numbers, 
the (closed tensor) dg-category corresponding to a $C^{\infty }$-manifold is  quasi-equivalent 
to the dg-category of flat bundles on the manifold ( \cite[section 3]{sim}). Here, the dg-category of flat bundles on a manifold $M$ is such that
\begin{itemize}
\item its objects are (finite rank) flat vector bundles $(V,D)$ on $M$, where  $V$ is a vector bundle on $M$ and $D$ is a flat connection on $V$, and
\item its complex of morphisms between two flat bundles $(V,D), (V',D')$ is the de Rham complex of $M$ with coefficients in the flat bundle $(\inhom (V,V'),D_{\inhom }\ )$, where $\inhom(V,V')$ is the hom-vector bundle between $V$ and $V'$ and $D_{\inhom}$ is the flat connection which is induced by $D,D'$.
\end{itemize}
This dg-category has natural closed tensor structure. If $M$ is simply connected, all flat bundles are trivial and this dg-category is essentially the same as the de Rham algebra. For rational coefficients, we construct the corresponding dg-category, using polynomial de Rham forms and finite rank rational local systems instead of flat bundles.  
\subsubsection{Main results}
Let $\dgcl$ be the category of small closed tensor dg-categories. The first main result is the following.
\begin{thm}[Theorem \ref{modelstrondgcl}]\label{result1}
The category $\dgcl$ has a model category structure where weak equivalences are quasi-equivalences. 
\end{thm}
We extract this result from  a theorem of Tabuada which states the category of small dg-categories admits a model category structure using Quillen's path-object argument (lifting argument, see Theorem \ref{pathobjarg}). The main problem is construction of free functor i.e., a left adjoint of the forgetful functor from the category of small closed tensor dg-categories to the category of small dg-categories. \\

For each pointed simplicial set $K$, we construct a closed tensor dg-category $\tdrk$ with a fiber functor. This construction is functorial in the  contravariant sense. Let $\dgclp$ be the category of small closed tensor dg-categories with fiber functors. The second main result is the following.
\begin{thm}[Theorem \ref{mainthm}]\label{result2}
Let $\ssetpfq$ be the category of connected pointed simplicial sets whose fundamental groups are finite and whose higher homotopy groups are uniquely divisible and finite dimensional as $\bb{Q}$-vector spaces. \\
\textup{(1)} There exists a full subcategory $\tannf$ of $\dgclp$ and the functor $K\mapsto \tdrk$ induces an equivalence of homotopy categories:
$\ho(\ssetpfq)\simeq\ho (\tannf)^{op}$.\\
\textup{(2)} Let $K$ be a simplicial set whose fundamental group is finite and whose higher homotopy groups are  Abelian groups of finite rank. The adjunction map 
\[
K\longrightarrow \bb{R}\mathit{Sp}_0\tdrk,
\]
where $\bb{R}\mathit{Sp}_0$ is a right adjoint of $\tdr$, is a fiberwise rationalization of $K$.
\end{thm}
The points of the proof of this result are as follows.
\begin{itemize}
\item For a finite group $G$, the unit morphism $K(G,1)\to\bb{R}\mathit{Sp}_0\tdr (K(G,1))$ is a weak equivalence of simplicial sets. 
\item Let $L$ be a simplicial set whose fundamental group is finite and whose higher homotopy groups are finite dimensional rational vector spaces, and  $\widetilde L\to L\to K(\pi _1(L),1)$ be a  homotopy fiber sequence where the map $L\to K(\pi _1(L),1)$ induces an isomorphism of $\pi_1$. The corresponding sequence $\tdr (K(\pi _1(L),1))\to \tdr (L)\to \tdr(\widetilde L)$ is a homotopy cofiber sequence in the category of closed tensor dg-categories with fiber functors. 
\end{itemize}
In the infinite fundamental group case, we cannot expect the rational homotopy types in the above sense are recovered from the corresponding closed tensor dg-categories. It is likely that these closed tensor dg-categories are equivalent to To\"en's schematic homotopy types (see \cite{champs, kpt} and subsection \ref{background}) but we don't discuss this in the present paper.

\subsubsection{Relation with equivariant differential graded algebras}\label{relation}
 We shall mention the relation between our formulation and the formulation using equivariant dg-algebras (see \cite{brownszc}).  An equivariant (commutative) dg-algebra is, by definition, a commutative dg-algebra with a group action. Let us consider the pointed case. Let $K$ be a  possibly non-simply connected  pointed simplicial set. We take the universal covering $\widetilde K\to K$. The corresponding polynomial de Rham algebra $A_{dR}(\widetilde K)$ has a natural action of $\pi_1(K)$ induced by the action of $\pi _1(K)$ on $\widetilde K$. Let $\widetilde{A_{dR}}(K)$ denote the equivariant dg-algebra $(\pi_1(K),A_{dR}(\widetilde K))$. Under the finiteness conditions on the higher homotopy groups, one can recover  the rational homotopy type of $K$. In the finite fundamental group case, the closed tensor dg-category $\tdrk$ and $\widetilde{A_{dR}}( K)$ are equivalent in the following sense.  The objects of $\tdr (K)$ are
 by definition, finite rank $\bb{Q}$-local systems  on $K$ or equivalently, finite dimensional  $\bb{Q}$-representations of $\pi _1(K)$. Let $\mathbf{1}$ be a trivial 1-dimensional representation and $V_r$ be the regular representation (see definitions below Lemma \ref{cofibrancy}). Consider the complex of morphisms $A:=\homo _{\tdrk}(\mathbf{1},V_r)$. The pointwise multiplication $V_r\otimes V_r\to V_r$ induces a structure of commutative dg-algebra on $A$, the right action of $\pi_1$ on $V_r$ induces an action on $A$ and one can see that $A$ is isomorphic to $\widetilde{A_{dR}}( K)$.  On the other hand, one can construct a closed tensor dg-category which is equivalent (in the categorical sense) to $\tdr (K)$ from $\widetilde{A_{dR}}( K)$.\\
\indent  We have the following diagram of categories.
\[
\xymatrix{
&&(\eqalgfp)^{op} \ar[d]^{\ms{T}}_{}="t"\\
\ssetpfq \ar[urr]^{\widetilde{A_{dR}}}\ar[rr]_{T_{dR}}^{}="s" &&(\tannf)^{op}\ar @/^/@{=>}"s";"t"^{\Phi}}
\]
Here, $\eqalgfp$ is the category of 1-connected augmented equivariant dg-algebras of finite types (see Definition \ref{defofeqalg}) and $\T$ is a functor. 
Comparison results are summarized as follows.
\begin{thm}[Theorem \ref{eqalg}, Proposition \ref{naturaltrans}]\label{result3} 
\  \\ 
\textup{(1)} There exists a natural transformation 
\[
\Phi:\tdr\Longrightarrow \ms{T}\circ \widetilde{A_{dR}}:\ssetpfq\longrightarrow(\dgclp)^{op}
\]  
such that for each $K\in\ssetpfq$, $\Phi _K$ is an equivalence of categories (which underlie closed tensor dg-categories).\\
\textup{(2)} The  functor $\T$ induces an equivalence of  homotopy categories $:\ho (\eqalgfp)\stackrel{\sim}{\longrightarrow}\ho(\tannf)$ 
\end{thm}
One can prove the functor $\widetilde{A_{dR}}$ induces an equivalence of homotopy categories independently of Theorem \ref{result2} and  \ref{result3}, though we don't prove in this paper. If we assume this equivalence, Theorem \ref{result2} follows from Theorem \ref{result3}. Our way of the proof of Theorem \ref{result2} is not the shortest one  but we take the way in order to understand our algebraic objects. In the proof of Theorem \ref{result3}, (2) we use internal hom functor.\\
\indent If the fundamental group is infinite, the closed tensor dg-category is not equivalent to  the equivariant  dg-algebra in any sense because we consider only finite rank local systems.

\subsubsection{Organization of the paper}
 	We review the contents of this paper. The main body is the second and third sections. In the second section, we prove the existence of a model category structure on the category of closed tensor dg-categories (Theorem \ref{result1}).  In \ref{preliminaries}, we give the definition of closed tensor dg-category and gather known results which is used in the proof.  In \ref{free} we  construct the free functor, which is necessary for the path object argument. In \ref{proofofmodel} we complete the proof. \\
\indent	 In the third section we prove the equivalence theorem \ref{result2}. In \ref{scdga}, we define a Quillen pair between the category of simplicial sets and the opposite category of the category of closed tensor dg-categories. Its left Quillen functor is the above $\tdr$.  In \ref{tannandeqalg} we compare closed tensor dg-categories with equivariant dg-algebras. We prove Theorem \ref{result3} and more rigid result (see Theorem \ref{eqalg}, (1), (2)). The main tool used here is Tannakian theory summarized in Theorem \ref{tdforfinite}. We provide some explicit examples of closed tensor dg-categories.  We also prove  a lemma about homotopy pushout. In \ref{suldr}  we prove Theorem\ref{result2}. \\
\indent One can read the third section independently of the second section if he or she assume Theorem \ref{modelstrondgcl}, Corollary \ref{modelstrondgclp}. \\
\indent  Arguments in this paper are all elementary except the language of model category theory.

\subsubsection{Background}\label{background} 
Our motivation is the application of non-Abelian Hodge theory to the topology of complex projective manifold. In simply connected case, combined with Sullivan's result, Hodge theory gives mixed Hodge structures on rational homotopy groups and rational minimal models  of compact K\"ahler manifolds and complex quasi-projective manifolds and then, the mixed Hodge structures give restrictions to the topology of them.(see \cite{suldel,mor}). As a generalization of these results to the non-simply connected case,  the application of non-Abelian Hodge theory is studied by Katzarkov, Pantev, To\"en \cite{toen} and Pridham \cite{prid0,prid}. Non-Abelian Hodge theory  states 
quasi-equivalence between the dg-category of flat bundles and the dg-category of semistable Higgs bundles 
with vanishing Chern numbers on  complex projective manifolds (see \cite[Section 3]{sim}).
In \cite{toen,prid} the authors define and construct  "mixed Hodge structure" on some algebraic object encoding homotopical data of a complex manifold, by using non-Abelian Hodge theory, then in \cite{toen},  restrictions to the homotopy types of complex projective manifolds are given and in \cite{prid}  mixed Hodge structures on real homotopy groups of them are constructed under some assumptions. In \cite{toen}, the algebraic objects are  schematic homotopy types, which are higher stacks  and in \cite{prid} the ones are the pro-algebraic homotopy types, which are simplicial affine group schemes (see also \cite{pridham}). Our algebraic objects can be  an alternative approach to these problems and we think the use of dg-category is more natural because the dg-category of flat bundles is a natural extension of the de 
Rham algebra and the dg-category naturally appears in non-Abelian Hodge theory.

\subsection{Notation and terminology} 
Throughout this paper, $k$ denotes a field of characteristic 0 and $\bb{Q}$ denotes 
the field of rational numbers. 
\subsubsection{dg-categories}

All complexes are defined over $k$ and have cohomological grading.
$\ms{C}(k)$ denotes the symmetric monoidal category of unbounded complexes. By a differential graded (dg-) category, we mean a category enriched over $\ms{C}(k)$ (See \cite{keller}). $\underline{\ms{C}}(k)$ denotes the dg-category of unbounded complexes i.e., the self-enrichment of $\ms{C}(k)$. For a 
(dg-)category $C$, $\ob (C)$ denotes the  set of objects of $C$ and for a (dg-)functor $F:C\to D$, $\ob (F):\ob (C)\to \ob (D)$ denotes the fuction given by $F$. If $C$ is a category, 
$\homo _C(c,c')$ stands for the set of morphisms between $c$ and $c'$. If $C$ is 
a dg-category, the same symbol denotes the complex of morphisms. 
We always identify $k$-linear categories with 
dg-categories concentrated in degree 0. 
Commutative dg-algebra is abbreviated to cdga and dg-category to dgc.  \\

We denote by $\dgcu$ the category of small dg-categories  and dg-functors  between them (see \cite{keller}), 
by $\dgc$ the full subcategory of $\dgcu$ consisting of dg-categories $C$ such that 
$\homo _C^n(c,c')=0$ for any $c,c'\in \ob (C)$ and for any $n< 0$.\\

 For $C\in \dgcu$ let $Z^0(C)$ (resp. 
$H^0(C)$) denote 
the category whose objects are those of $C$ and whose sets of morphisms consist of 0-th cocycles 
(resp. 0-th cohomology classes) of complexes of morphims of 
$C$.  
We regard  $Z^0$ and $H^0$ as functors
\[
Z^0, H^0:\dgcu\longrightarrow \ms{Cat},
\]
where $\ms{Cat}$ is the category of small categories. It is clear what $Z^0$ and $H^0$ assign to each morphism of $\dgcu$.
A morphism in $Z^0(C)$ is 
said to be a chain morphism in $C$. 
A morphism in $C$ is said to be an isomorphism if  it is a chain morphism and has an inverse 
which is also a chain morphism. 
$\mor (C)$ stands for the set of all homogeneous morphisms of $C$, i.e., $\mor (C):=\bigsqcup\limits_{(c,c')}
\bigsqcup\limits_n\homo _C^n(c,c')$. \\

Let $F$, $G:C\to D$ be two dg-functors. A natural transformation 
$\alpha :F\Rightarrow G$ requires that for each $c\in \ob (C)$ 
$\alpha _c$ is a chain morphism and compatible with all morphisms of 
all degrees. \\

An equivalence (resp. a quasi-equivalence)  
between dg-categories is a dg-functor which induces an equivalence between $Z^0$'s (resp. $H^0$'s) 
and isomorphisms (resp. quasi-isomorphisms) of 
the complexes of morphisms.\\

\indent Let $C,D\in\dgc$. 
\begin{itemize}
\item $C\boxtimes D$ denotes 
a dgc defined as follows.
\begin{itemize}
\item $Ob(C\bt D)=\{(c,d)\, |\, c\in C,\ d\in D\}$ and

\item $Hom_{C\bt D}((c_0,d_0),(c_1,d_1))=Hom_C(c_0,c_1)\ot_k Hom_D(d_0,d_1)$
with the composition given by $(f'\ot g')\circ (f\ot g)=(-1)^{\deg g'\cdot \deg f}
(f'\circ f)\ot (g'\circ g)$.
\end{itemize}
\item $\mca{T}_{C,D} :C\bt D\rightarrow D\bt C$ denotes the morphism defined by 
$(c,d)\mapsto (d,c)$ and $f\ot g\mapsto (-1)^{\deg f\cdot \deg g}g\ot f$. 
\item $C^{op}$ denotes the opposite dg-category of $C$ whose composition is defined by 
$g\circ f:=(-1)^{\deg g\cdot\deg f}f\circ g $, where the composition of right hand side 
is that in $C$. 
\item We define a dg-functor $\homo _C(-,-):C^{op}\bt C\to \underline{\ms{C}}(k)$ by 
$\homo _C(c, c')=\homo _C(c,c')$ for $(c, c')\in \ob (C^{op}\bt C)$ and 
\[
\homo_C (f\otimes g)(\alpha )=(-1)^{\deg f(\deg g+\deg \alpha )}g\circ\alpha\circ f
\]
for $f\otimes g\in \mor (C^{op}\bt C)$.

\item $C\times D$ denotes the product in the category $\dgc$. Explicitly, $\ob(C\times D)=\ob(C)\times \ob(D)$, $\homo_{C\times D}((c_0,d_0),(c_1,d_1))=\homo_C(c_0,c_1)\oplus\homo_D(d_0,d_1)$. 
\end{itemize}
Clearly these constructions are functorial.
\vspace{1\baselineskip}

\indent We denote by $\dggr$ the category of non-negatively graded dg-graphs. Its objects are 
directed graphs whose edges have  structures of non-negatively graded complexes and its morphisms are 
morphisms of directed graphs which induce chain maps on edges.  
$\fcat :\dggr \to \dgc $ denotes the free functor of \cite[Section 5]{tab3}. Explicitly,  
$\ob (\fcat (G))=\vtx (G)$, the set of vertices of $G$, and 
\[
\homo _{\fcat (G)}(v,v'):=\left\{
\begin{array}{ll}
k\!\cdot\! id _v\oplus \bigoplus \limits_{ v_1,\dots ,v_l,l\geq 0}\ed (v_l,v')\ot _k\cdots\ot  _k\ed (v,v_1)&\textrm{if }v=v'\\
\bigoplus \limits_{v_1,\dots ,v_l\geq 0}\ed (v_l,v')\ot_k\cdots\ot  _k\ed (v,v_1)&\textrm{otherwise} \\
\end{array}\right. 
\]
where $\ed$ stands for the complex of edges.  
$\fcat $ is a left adjoint of the forgetful functor $\dgc\to \dggr$.\\
 
 A non-unital dg-category is a dg-graph $G$ which is equipped with associative composition $\ed(v',v'')\otimes _k\ed(v,v')\to\ed (v,v'')$ for each $v,v',v''\in \vtx (G)$ (but without units). An ideal of a dg-category $C$ is a non-unital dg-subcategory $I$ of $C$ such that it 
contains all objects of $C$ and if $f$ and $g$ are morphisms of $C$, one of them is in 
$I$, and $g\circ f$ exists, $g\circ f$ is in $I$. 
If $I$ is an ideal of $C$, then a dg-category $C/I$ is defined by 
$\ob (C/I)=\ob (C)$ and $\homo _{C/I}(c,c')=\homo _C(c,c')/\homo _I(c,c')$.  
For  a subset $S$ of $\mor (C)$, the ideal generated by $S$ is the 
smallest ideal which contains $S$.\\

\subsubsection{Model categories}\label{modelcategories}

Our notion of model category is that of \cite{hov}.  Let $\ms{M}$ be a model category. \\

\indent Let $\emptyset\in \ms{M}$ be an initial objects. The over category $\ms{M}/\emptyset$ has a model category structure where equivalences, fibrations and cofibrations are detected by underlying morphisms of $\ms{M}$ (see \cite[Remark 1.1.7, Proposition 1.1.8]{hov}). We always  regard the category $\ms{M}/\emptyset$ as a model  category by this structure.\\

\indent The notions of path object and right homotopy are found in \cite[DEFINITION 1.2.4]{hov}.\\

\indent The notion of homotopy pushout squares in $\ms{M}$ is found in \cite[p.184]{hov} and if $\ms{M}$ is a pointed category, a sequence $X\stackrel{f}{\to}Y\stackrel{g}{\to}Z$ of morphisms of $\ms{M}$ with $g\circ f=0$ is said to be a homotopy cofiber sequence if the commutative square
\[
\begin{CD}
X@>f>>Y\\
@VVV@VgVV\\
*@>>>Z\\
\end{CD}
\]
is a homotopy pushout square. Here $*$ is a terminal (and initial) object of  $\ms{M}$. The notion of homotopy fiber sequences is its dual notion.\\

\indent $\ho (\ms{M})$ denotes 
the homotopy category of $\ms{M}$. Let $\ms{M}'$ be a full subcategory of $\ms{M}$ stable under weak equivalences of 
$\ms{M}$. We denote by $\ho (\ms{M}')$ the full subcategory of $\ho (\ms{M})$ consisting of objects of $\ms{M}'$. It is easy to see $\ho (\ms{M}')$ is (isomorphic to) the localization of $\ms{M}'$ obtained by inverting weak equivalences in $\ms{M}'$ \\

\indent $\sset$ 
stands for the category of simplicial sets. We regard $\sset$ as a model category with the usual model structure (see \cite{hov}).  $\ssetp$ stands for the category of pointed simplicial sets. $\ssetpf$ (resp. $\ssetpfq$) denotes the full subcategory of $\ssetp$ consisting of connected $K$ with $\pi_1(K)$ finite, $\pi _n(K)$ being an Abelian group of finite rank (resp. with $\pi_1(K)$ finite, $\pi _n(K)$ uniquely divisible and finite dimensional as a $\bb{Q}$-vector space) for each $n\geq 2$.  Let 
$K\in\sset$. $\Delta (K)$ denotes the category of simplices of $K$ of \cite[Chapter 3]{hov}:
\begin{itemize}
\item An object of $\Delta (K)$ is a simplex of $K$, i.e., $\ob (\Delta (K))=\bigsqcup\limits_{n\geq 0}K_n$,
\item for $\sigma\in K_n$, $\tau\in K_m$ a morphism $a:\sigma\to\tau $ is a morphism $a:[n]\to[m]\in \Delta$ 
such that $a^{*}(\tau )=\sigma $
\end{itemize}
where $\Delta$ is the category with objects $[l]=\{0,\dots,l\}$ for $l\geq 0$, 
and weakly order-preserving maps.

\section{ Model of homotopy theory of closed tensor dg-categories}
\indent The purpose of this section is to give a model structure on the category of small closed tensor 
dg-categories (see Theorem \ref{modelstrondgcl}). This result is  a foundation of arguments of next section.

\subsection{Preliminaries}\label{preliminaries}

\subsubsection{Closed tensor dg-categories}
The following definition is a variant of the usual definition of closed symmetric monoidal category (see \cite[Definition 4.1.12]{hov} or \cite[P.180]{mac}, where the author call it closed category) in the differential graded context. 
\begin{defi}\label{defofcldgc}
\textup{(1)} Let $C$ be an object of $\dgc$. A \textup{closed tensor structure} on $C$ 
is a 11-tuple
\[
((-\ot -),\, \mb{1},\, a,\, \tau ,\, u,\, \inhom,\,  \phi,\,  (-\oplus -),\, s_1,\, s_2,\, \mb{0})
\]
consisting of 
\begin{itemize}
\item a morphism $(-\ot -):C\bt C\longrightarrow C \in \dgc$,
\item a distinguished object $\mb{1}\in C$,

\item natural isomorphisms 
\begin{align*}
a :((-\ot -)\ot -)&\Longrightarrow  
(-\ot (-\ot -))\,\, \, :(C\bt C)\bt C\cong C\bt (C\bt C)\longrightarrow C,  \\
\tau:(-\ot -)&\Longrightarrow (-\ot -)\circ \mca{T}_{C,C}:C\bt C\longrightarrow C,\\
u:(-\ot \mb{1})&\Longrightarrow id_C\qquad\qquad\ \,:C\longrightarrow C
\end{align*}
satisfying
usual coherence conditions on associativity, commutativity and unity, see \cite[pp.251]{mac},
\item a morphism $\inhom :C^{op}\bt C\longrightarrow C \in \dgc$,
\item a natural isomorphism 
\begin{align*}
\phi :\homo _C(-\ot -,-)\Longrightarrow \homo _C(-&,\inhom (-,-)):C^{op}\bt C^{op}\bt C\longrightarrow \ms{C}(k),
\end{align*}
\item a morphism $(-\oplus -):C\times C\longrightarrow C\in \dgc$,
\item two natural transformations 
\[
\xymatrix{P_1\ar@{=>}[r]^{s_1\ \ \ }&(-\oplus -)&
P_2:C\times C\ar@{=>}[l]_{s_2}\ar[r]& C,}
\] where $P_i:C\times C\longrightarrow C$ is 
the $i$-th projection, such that the induced morphism $s_1^*\times s_2^*:\homo _C(c_0\oplus c_1,c')
\longrightarrow \homo_C(c_0,c')\times \homo_C(c_1,c')$ is an isomorphism (i.e., $c_0\oplus c_1$ is a coproduct), and
\item a distinguished object $\mb{0}\in \ob (C)$ such that $\homo _C(\mb{0},c)=0$ for any 
$c\in \ob (C)$. 
\end{itemize}
We call $(-\otimes -)$ a \textup{tensor functor} and $\inhom$ a \textup{internal hom functor}.\\
\textup{(2)} A \textup{closed tensor dg-category} is an object $C$ of $\dgc$ equipped with a closed tensor structure. For two closed tensor dg-categories $C,D$, a \textup{morphism of closed tensor dg-categories} is a morphism $F:C\rightarrow D$ 
of dgc's which preserves all of the above structures. For example,
$F(c\otimes d)=F(c)\otimes F(d)$ (not only naturally isomorphic), $F(\tau_{c,c'})=\tau_{Fc,Fc'}$ 
and $F(\mb{1})=\mb{1}$. 
We denote by $\dgcl$ the category of closed tensor dgc's. A \textup{closed tensor category} is a closed tensor dg-category whose complexes of morphisms are concentrated in degree 0 and a \textup{morphism of closed tensor  categories} is the same as a morphism of closed tensor dg-categories.
\end{defi}
We apply the notions of equivalence and quasi-equivalence to 
objects of $\dgcl$ via the forgetful functor $\dgcl\to\dgc$. Note that an equivalence  in $\dgcl$ doesn't always have a  
quasi-inverse which is a morphism of $\dgcl$. We say two objects of $\dgcl$ are equivalent if  
they can be connected by a finite chain of equivalences in $\dgcl$.\\

\indent Let $T:C\bt C$ (resp. $C^{op}\bt C$) $\rightarrow C$ be a dg-functor. 
A $T$-closed ideal is an ideal $I$ of $C$ such that $T(f,g)$ is in $I$ if one of $f$ and $g$ is in $I$. 
If $I$ is $T$-closed ideal, $T$ induces a functor $\oli{T}:(C/I)\bt (C/I)$ (resp. $(C/I)^{op}\bt (C/I)$) 
$\to C/I$. There
is an obvious notion of $T$-closed ideal generated by $S$. 
\subsubsection{A model category structure on $\dgcu $}

\indent The following theorem is the main result  of \cite{tab1} which is crucial for our argument.
\begin{thm}[\cite{tab1}]\label{dgcu}
The category $\dgcu$ has a cofibrantly generated model structure where 
weak equivalences and fibrations are
defined as follows.
\begin{itemize}
\item A morphism $F:C\rightarrow D\in\dgcu$ is a weak equivalence if and only if
it is a quasi-equivalence.
\item A morphism $F:C\rightarrow D\in\dgcu$ is a fibration if and only if it satisfies 
the following two conditions.
\begin{itemize}
\item For $c,c'\in \ob (C)$ the morphism $F_{(c,c')}:\homo _C(c,c')\to \homo _D(Fc,Fc')$ is an epimorphism. 

\item For any $c\in \ob (C)$ and any isomorphism $f:Fc\to d'\in H^0(D)$, there exists an isomorphism 
$g:c\to c' \in H^0(C)$ such that $H^0(F)(g)=f$.
\end{itemize}
\end{itemize}
\end{thm}
\subsubsection{Path object argument}
The path object argument is due to Quillen \cite{homotopicalalg} and the following form is found in \cite[section 5]{tab2}. 
\begin{thm}[Path object argument, \cite{homotopicalalg,tab2}]\label{pathobjarg} 
Let $\ms{M}$ be a category with all small limits and colimits and $\ms{N}$ be a cofibrantly generated model category. Let $\mca{U}:\ms{M}\longrightarrow \ms{N}$ be a functor which commutes with all filtered colimits.
Suppose that
\begin{itemize}
\item all objects of $\ms{N}$ are fibrant,
\item $\mca{U}$ possesses a left adjoint $\mca{L}:\ms{N}\longrightarrow \ms{M}$, and 
\item There exists an endo-functor $\mca{P}:\ms{M}\longrightarrow \ms{M}$ and a sequence of natural transformations 
\[
\xymatrix{\mathit{Id}_{\ms{M}} \ar@{=>}[r]^i
& \mca{P}\ar@{=>}[r]^{d_0\times d_1\ \ }&\Pi :\ms{M}\ar[r] &\ms{M},}
\] 
where $\Pi$ is defined by $\Pi(X)=X\times X$, such that for each object $X\in \ms{M}$ the composition $(d_0\times d_1)_X\circ i_X$ is the diagonal $X\to X\times X$ and the diagram $\mca{U}X\stackrel{\mca{U}i}{\to}\mca{UP}X\stackrel{\mca{U}d_0\times \mca{U}d_1}{\to}\mca{U}\Pi X=\mca{U}X\times \mca{U}X$ is a path object in the sense of \cite[Definition 1.2.4]{hov}.
\end{itemize}
Then, there exists a cofibrantly generated model category structure on $\ms{M}$ such that a morphism $f:X\to Y\in \ms{M}$ is a weak equivalence (resp. a fibration) if and only if  $\mca{U}f:\mca{U}X\to \mca{U}Y$ is a weak equivalence (resp. a fibration) in $\ms{N}$. 
\end{thm}
We call a functor $\mca{P}$ satisfying the above condition a functorial path object.

\subsection{Free construction}\label{free}
\indent In order to use path object argument, we shall construct a left adjoint  
\[
\fcl ^{\geq 0}  :\dgcu\longrightarrow \dgcl .
\]
of the forgetful functor 
$\mca{U}:\dgcl \to \dgcu$. We call this functor the free functor. (We use the definite article in the sense that it is unique up to natural isomorphisms.)  We divide construction of $\fcl^{\geq 0}$ into construction of two functors, 
\[
\mca{T}^{\geq 0}:\dgcu\longrightarrow \dgc \ \ \  \textrm{and}\ \ \  \fcl :\dgc\longrightarrow \dgcl .
\] 
Here, $\mca{T}^{\geq 0}$ is  a left adjoint of the inclusion functor 
$\mca{I}:\dgc\to\dgcu $ and $\fcl$ is a left adjoint of the forgetful functor: $\dgcl\to\dgc$, which we call the free functor too. If these two functors exist, it is  clear that $\fcl^{\geq 0}=\fcl\circ\mca{T}^{\geq 0}$.\\ 
\indent We first define $\mca{T}^{\geq 0}$. Let $C\in \dgcu$. For each $c, c'\in \ob(C)$, let $M(c,c')\subset \homo_C(c,c')$ be the subcomplex generated by homogeneous elements of the form 
\[
\sum_if_{i,1}\circ\cdots\circ f_{i,k_i},  
\]
where each $f_{i,j}$ is a homogeneous morphism and for each $i$, at least one of $f_{i,j}$'s has negative degree.
We put 
\[
\ob(\mca{T}^{\geq 0}C)=\ob(C),\ \ \ \ \homo_{\mca{T}^{\geq 0}C}(c,c')=\homo_C(c,c')/M(c,c')
\]
and define the composition of $\mca{T}^{\geq 0}C$ from that of $C$. One can easily see the construction $C\longmapsto \mca{T}^{\geq 0}C$ is functorial and satisfies the required property.\\

\indent The idea to construct $\fcl :\dgc\to\dgcl$ is also elementary: Attaching necessary objects and morphisms and 
taking quotients by necessary relations. To make this precise, we first define the following: 
\subsubsection{Universal dg-categories}
\indent In this sub-subsection, 
 we define a dgc which is initial among 
the dgc's having given objects, given morphisms, given relations, and a morphism from given dgc. 
Suppose the following data are given:
\begin{itemize}
\item  a set $\sob$,

\item a set $\smor$ of non-negatively graded complexes (i.e., $\smor$ is a subset of $\ob(\ms{C}(k))$).

\item two functions $s,\, t:\smor\to\sob $ which we call the source function and target function, 
respectively,  

\item a dgc $C$ and a function $\mf{o}:\ob (C)\to\sob $.
\end{itemize}
Let $\bigsqcup \smor$ be the set of homogeneous elements of complexes belonging to $\smor$, i.e., 
$\bigsqcup\smor =\bigsqcup\limits_{n\geq 0}\bigsqcup\limits_{H\in\smor}H^n$. \\
\indent Let $\mor (C,\smor )$ be the set of "formal morphisms generated by $\mor (C)$ and $\bigsqcup\smor$". More precisely, an element of the set $\mor (C,\smor )$ is a formal linear combination of formal compositions:
\[
\sum _n c_n(\alpha _{n ,k_n}\circ\cdots\circ \alpha _{n ,1})
\]
such that  $c_n\in k$ and each $\alpha _{n,i}$ is an element of the union 
\[
\mor (C)\sqcup \bigsqcup \smor \sqcup \{ id_x\}_{x\in\sob}
\]
($id_x$ is a formal symbol),  $t(\alpha _{n,i})=s(\alpha _{n,i+1})$ for each $n$ and each $i=1,\dots k_n-1$, and $t(\alpha _{n, k_n})=t(\alpha _{n', k_{n'}})$, $s(\alpha _{n, 1})=s(\alpha _{n',1})$ for each $n,n'$. Here, if $\alpha$ is 
a morphism in $C$, $s(\alpha)$ (resp. $t(\alpha)$ ) denotes the image of the source of $\alpha$ 
(resp. the target of $\alpha$) by $\mf{o}$.\\
\indent   Suppose the following additional datum 
is given:
\begin{itemize}
\item $\srela $, a subset of $\mor (C,\smor)$
\end{itemize}
\begin{defi}
With the above notations, the \textup{universal dgc associated to} $(C,\sob,\smor,\srela,s,t,\mf{o})$ is a 5-tuple 
\[
(D,\,I_D,\,\mf{o}_D,\, \mf{m}_D,\,\{f_{D,H}\}_{H\in \smor})
\]
consisting of  
\begin{itemize}
\item a dgc $D$, written $C[\sob , \smor]/\srela $, 
\item a morphism of dgc's $I_D:C\to D$,
\item a function $\mf{o}_D:\sob\to \ob (D)$,  
\item a function $\mf{m}_D:\smor \to \underline{\mor}(D )$ 
where 
$\underline{\mor}(D)=\{\,\homo _{D}(x,y)\,|\,(x,y)\in\ob (D )^{\times 2} \,\}$, and
\item a family of homomorphisms of complexes $\{f_{D,H}:H\to\mf{m}_D(H)\}_{H\in \smor}$
\end{itemize}
satisfying the following conditions.
\begin{itemize}
\item It must satisfy obvious consistency conditions. Firstly, the diagram 
\[
\begin{CD}
\smor @>\mf{m}_D>>\underline{\mor }(D)\\
@Vs \textrm{ (resp. $t$)}VV@VVs_D\textrm{ (resp. $t_D$)}V\\
\sob @>\mf{o}_D>>\ob (D),\\
\end{CD}
\]
where the right  vertical arrow is the source function
 (resp. the target function) of $D$, commutes. 
\item  Secondly, the composition 
$\ob (C)\stackrel{\mf{o}}{\to}\sob\stackrel{\mf{o}_D}{\to} \ob (D)$ is equal to the function  $\ob(I_D):\ob(C)\to \ob(D)$.

\item The function $\mor(C, \smor)\rightarrow \mor (D)$ defined from $I_D$, $\mf{m}_D$ and $f_{D,H}$'s 
takes $\srela$ to zeros.

\item It has a universal property. If a 5-tuple 
\[
(E,\,I_E:C\to E,\,\mf{o}_E:\sob\to \ob (E),\, \mf{m}_E :\smor \to \underline{\mor}(E),\,\{f_{E,H}:H\to\mf{m}_E(H)\})
\]
which satisfies all of the above conditions where $D$ is replaced with $E$ is given, 
there exists a unique morphism of dgc's $D\to E$ which preserves all of the above structures.
\end{itemize}
\end{defi}

\indent We shall construct the universal dgc.
If we have constructed such a dgc for the case $\srela=\emptyset$, then for general case, we only have  to put $C[\sob , \smor]/\srela =D_0/I$ where $D_0=C[\sob , \smor]/\emptyset$ and $I$ is the 
ideal of $D_0$ generated by  the image of $\srela$ by the function $\mor (C,\smor)\rightarrow \mor (D_0)$. 
So we may assume  $\srela=\emptyset$.\\ 
\indent We first define a dg-graph $A$ by 
\begin{itemize}
\item $\vtx (A)=\sob$, 
\item $\ed (x,y)=\bigoplus\limits_{
\begin{smallmatrix}
H\in\smor,\\
s(H)=x,\ t(H)=y
\end{smallmatrix}
}H\ \ \ \oplus
\bigoplus\limits_{\mf{o}(c)=x,\mf{o}(c')=y}\homo _C(c,c')
$.
\end{itemize}
We consider the dgc $D':=\fcat (A)$. Let $[f]$ denote the morphism of 
$\fcat(A)$ corresponding to an edge $f\in A$. Let $J$ be the ideal of $D'$ generated by 
\[
R=\{\:[id_c]-id_{\mf{o}(c)},\,[g\circ f]-[g]\circ [f]\: |
\: c\in\ob (C),\,f,\,g\in\mor (C)\:\}.
\]
Put $D=D'/J$. The set of relations $R$ ensures that 
one can define a dg-functor $C\to D$ by $c\mapsto \mf{o}(c)$ and $f\mapsto [f]$.
The other data are defined obviously and it is clear that 
$D$ is the required universal dgc.

\subsubsection{Free closed tensor dg-categories}\label{freecl}
We shall construct the free functor
\[
\fcl:\dgc\longrightarrow \dgcl ,
\] 
i.e., a left adjoint of the forgetful functor. \\

\indent Let $S$ be a set and $\wcl (S)$ be the set of words generated by $S$ and  formal 
symbols $\mb{1}$, $\mb{0}$ with operations $\otimes$, $\inhom$ and $\oplus$. 
More precisely, $\wcl (S)$ is defined inductively as follows. Set
\begin{align*}
\wcl ^0(S)&=S\sqcup \{\mb{1},\mb{0}\}, \\  
\wcl ^p(S)&=\{\: x\ot y,\,\inhom (x,y), \,x\oplus y  \: |\:  x,\,y\in \wcl^{p-1}(S)\:\}\,\sqcup \,\wcl ^0(S).
\end{align*}
Let 
$i_0:\wcl ^0(S)\rightarrow \wcl ^1(S)$ be the inclusion and 
$i_p:\wcl ^p(S)\rightarrow \wcl ^{p+1}(S)$ be the map defined inductively 
by $i_p(x\otimes y)=i_{p-1}x\otimes i_{p-1}y$, $i_p(\inhom (x,y))=\inhom 
(i_{p-1}x,i_{p-1}y)$ and $i_p(x\oplus y)=i_{p-1}x\oplus i_{p-1}y$. We identify $\wcl^p(S)$ with a subset of $\wcl^{p+1}(S)$ via $i_p$ and set 
\[
\wcl (S):=\bigcup_{p\geq 0}\wcl^p(S).
\]   Note that  for example, the operation $\otimes$ 
doesn't satisfy the associativity low and $\mb{1}$ doesn't play any special 
role yet. \\

Let $C\in\dgc$. To construct the free  closed tensor dg-category $\fcl (C)$ associated to $C$, 
we first construct the following \textbf{data}:
\begin{itemize}
\item a sequence of dgc's 
$C=D_{-1}\stackrel{I_{-1}}{\to}\cdots\stackrel{I_{j-1}}{\to}D_j\stackrel{I_{j}}{\to}\cdots$ such that $\ob (D_j)=\wcl(\ob(C))$ for $j\geq 0$, 

\item a morphism of dgc's $T_j:D_j\bt D_j\rightarrow D_{j+1}$ for each $j\geq -1$,

\item isomorphisms $a_{x,y,z}:(x\otimes y)\otimes z\rightarrow
x\otimes (y\otimes z)$, $\tau_{x,y}:x\otimes y\rightarrow y\otimes x\in \mor (D_0)$ whose inverse is 
$\tau_{y,x}$, and 
$u_x: x\otimes \mb{1}\to x\in \mor (D_0)$
 for each $x,\, y$ and $z\in \ob (D_0)$

\item a morphism of dgc's 
$\inh_j:D_j^{op}\bt D_j\rightarrow D_{j+1}$ for each $j\geq -1$ and 

\item chain morphisms $\ev _x^y:\inhom (x,y)\otimes x\rightarrow y$ and 
$\coev _x^y:x\rightarrow \inhom (y,x\otimes y)$ in $D_0$, which we call an evaluation 
and a coevaluation, for each $x$, $y\in D_0$,

\item chain morphisms $s_i^{x_0, x_1}:x_i\to x_0\oplus x_1$, $p_i^{x_0,x_1}:x_0\oplus x_1\to x_i$ for $x_i\in\ob (D_0)$, for $i=0,1$ 
\end{itemize}
which satisfy the following \textbf{conditions}:

\begin{itemize}
\item $\ob (I_{-1})$ is the inclusion $\ob (C)\to \wcl (\ob (C))$ and for $j\geq 0$ $\ob (I_j)$ is 
the identity. 
\item For $j\geq -1$, $\ob (T_j)=\otimes $ and $\ob (\inh _j)=\inhom$ where the right hand side of 
each equation is the operation of $\wcl(\ob (C))$.

\item For $j\geq -1$, $I_{j+1}\circ T_j=T_{j+1}\circ (I_j\bt I_j):D_j\bt D_j
\rightarrow D_{j+2}$.

\item For $j\geq -1$, $I_{j+1}\circ \inh _j=\inh _{j+1}\circ (I_j^{op}\bt I_j):D_j^{op}\bt D_j
\rightarrow D_j$.
\item For $j\geq -1$, $(I_j,T_j,\inh _j)$ has a \textit{universal property} as follows. Suppose 
the following data are given.
\begin{itemize}
\item a closed tensor dgc $E$, 
\item a morphism of dgc's $I':D_j\to E$.

\end{itemize}
Suppose also these data satisfy the following conditions.
\begin{itemize}
\item $I'\circ T_{j-1}=(-\otimes_E-)\circ (I'\bt I')\circ (I_{j-1}\bt I_{j-1}):D_{j-1}\bt D_{j-1}\longrightarrow E$.
\item $I'\circ \inh _{j-1}=\inhom_E \circ (I'^{op}\bt I')\circ (I_{j-1}^{op}\bt I_{j-1}):(D_{j-1})^{op}\bt D_{j-1}\longrightarrow E$.
\item $I'$ takes  $a$'s, $\tau$'s, $u$'s, $\ev$'s, $\coev$'s, $s_i$'s, $p_i$'s $\mb{1}$, and $\mb{0}$ to the corresponding morphisms and objects.   
\end{itemize}
(If $j=-1$, these conditions are ignored.) 
Then if $j\geq 0$ (resp. $j=-1$), 
there is a unique morphism $\tilde I':D_{j+1}\to E$ such that $\tilde I'\circ I_j=I'$, 
$\tilde I'\circ T_j=(-\otimes_E-)\circ (I'\bt I')$ and $\tilde I'\circ \inh _j=\inhom_E\circ (I'^{op}\boxtimes I')$ (resp.
$\tilde I'$ takes  $a$'s, $\tau$'s, $u$'s, $\ev$'s, $\coev$'s, $s_i$'s, $p_i$'s $\mb{1}$, and $\mb{0}$ to the corresponding morphisms and objects ).
\end{itemize} 
\textbf{Construction}: The construction of the above data proceeds in induction. Suppose we have constructed the stage $p\geq -1$. 
In other words, we have constructed 
\begin{itemize}
\item a sequence of dgc's $D_{-1}\stackrel{I_{-1}}{\to}\cdots\stackrel{I_{p-1}}{\to}D_p$, 

\item a morphism of dgc's $T_j:D_j\bt D_j\rightarrow D_{j+1}$ for each $j\leq p-1$ and 

\item a morphism of dgc's $\inh_j:D_j^{op}\bt D_j\rightarrow D_{j+1}$ for each $j\leq p-1$ 
\end{itemize}
which satisfy the above conditions. 
We shall construct  $D_{p+1}$. The idea is to define $\sob$, 
$\smor $ and $\srela$ appropriately, and put $D_{p+1}=D_p[\sob, \smor]/\srela $.\\

Set $\sob := \wcl (\ob (C))$ and $\mf{o}:=id_{\wcl(\ob(C))}$. We shall define $\smor$. \\
\textit{Notation}: In the following, 
for a set of complexes $M$ (which will be a subset of $\smor$) we use the expression 
\[
M=\{\: M(x,y)\: |\: (x,y)\in M_{ob}\: \}.
\]
This means $M$ consists of $M(x,y)$'s, $M(x,y)$ is a complex whose source (resp. target) is 
$x$ (resp. $y$) and $(x,y)$ runs through $M_{ob}$, which is a subset of $(\sob )^{\times 2}$.\\

When $p=-1$, We define sets of complexes 
\[
\tens,\  \ass,\ \ass^{-1},\  \comm, 
\uni,\ \uni^{-1},\ \mathit{Int},\ \mathit{Ev},\ \mathit{Coev},\ \mathit{Inc}_i,\ \mathit{Proj}_i\ \ (i=1,2).
\]
Put
\begin{eqnarray*}
\tens (x_1\ot x_2,y_1\ot y_2)&=&\homo _{D_p}(x_1,y_1)\ot_k\homo _{D_p}(x_1,y_1),\\
\Int (\,\inhom (x_1, x_2),\,\inhom (y_1,y_2)\,)&=&\homo _{D_p}(y_1,x_1)\ot_k \homo _{D_p}(x_2,y_2)\\
\\
\ass (\,(x_1\ot x_2)\ot x_3,\, x_1\ot (x_2\ot x_3)\,)&=&k\cdot a_{x_1,x_2,x_3}, \\
\comm (x_1\ot x_2, x_2\ot x_1)&=&k\cdot \tau _{x_1,x_2}, \\
\uni (x\ot \mb{1},x)&=&k\cdot u_x, \\
\Ev (\inhom (x,y)\ot x, y)&=&k\cdot \ev _x^y, \\
\Coev (x,\inhom (y,x\ot y))&=&k\cdot \coev _x^y ,\\
\mathit{Inc}_i(x_i,x_0\oplus x_1)&=&k\cdot s_i^{x_0,x_1},\\
\mathit{Proj}_i(x_0\oplus x_1,x_i)&=&k\cdot p_i^{x_0, x_1}. 
\end{eqnarray*}
It will be clear that what $\tens_{ob},\dots ,\Coev _{ob}$ are. For example,
\[
\tens _{ob}=\{\: (x_1\ot x_2,y_1\ot y_2)\: |\: x_1,\,x_2,\,y_1,\,y_2\in\sob \: \}.
\]
We agree that 
$\ass^{-1}$ (resp. $\uni^{-1}$) is a copy of $\ass $ (resp. $\uni$) 
whose source and target functions are replaced with each other and we denote by $a'_{x_1,x_2,x_3}$ 
(resp. $u'_x$)  
the element of the complex belonging to $\ass^{-1}$ 
(resp. $\uni ^{-1}$) which corresponds to $a_{x_1,x_2,x_3}$ (resp. $u_x$). We set 
\begin{align*}
\smor = &\tens\sqcup \ass\sqcup \ass^{-1}\sqcup  \comm\sqcup 
\uni\sqcup \uni^{-1}\sqcup \mathit{Int}\sqcup \mathit{Ev}\sqcup \mathit{Coev}\\
&\sqcup \mathit{Inc}_1\sqcup \mathit{Inc}_2\sqcup\mathit{Proj}_1\sqcup \mathit{Proj}_2.
\end{align*}

To define $\srela $, we define three subsets of $\mor (D_p,\smor)$, $R_1$, $R_2$ and $R_3$ as follows. 
In the following, $T'_p (f_1,f_2)$ (resp. $\inh' _p(f_1,f_2)$) denotes the element 
$f_1\ot f_2\in \tens (x_1\ot x_2,y_1\ot y_2)$ (resp. $\in\Int (\,\inhom (x_1, x_2),\,\inhom (y_1,y_2)\,)$ ).

\begin{eqnarray*}
R_1&=&\left\{\:
\begin{array}{c}
 T'_p(g_1, g_2)\circ T'_p(f_1, f_2)-(-1)^{\deg g_2\cdot \deg f_1}
T'_p(g_1\circ f_1,g_2\circ f_2),\\
\\
id_{x_1\otimes x_2}-T'_p(id_{x_1}, id_{x_2})
\end{array}\:\right\},\\
\\
\\
R_2&=&\left\{\:
\begin{array}{c}
 \inh' _p(g_1, g_2)\circ \inh' _p(f_1, f_2)-(-1)^{(\deg g_1+\deg g_2)\deg f_1}
\inh' _p(f_1\circ g_1, g_2\circ f_2),\\
\\
id_{\inhom (x_1, x_2)}-\inh' _p(id_{x_1}, id_{x_2})
\end{array}\:\right\},\\
\\
\\
R_3&=&\left\{\:
\begin{array}{c}
a'_{x_1,x_2,x_3}\circ a_{x_1,x_2,x_3}-id_{(x_1\otimes x_2)\otimes x_3},\ 
a_{x_1,x_2,x_3}\circ a'_{x_1,x_2,x_3}-id_{x_1\otimes (x_2\otimes x_3)},\\
\\
u'_x\circ u_x-id_{x\otimes\mb{1}},\ u_x\circ u'_x-id_x,\ \ 
\tau _{x_2,x_1}\circ \tau _{x_1,x_2}-id_{x_1\otimes x_2}\\
\end{array}\:\right\}.
\end{eqnarray*}
\\
We set $\srela =R_1\sqcup R_2\sqcup R_3$ and  $D_{p+1}:=D_p[\sob ,\smor ]/\srela$. 
$I_p:D_p\to D_{p+1}$ is the structure morphism of the universal dgc. 
$R_1$ (resp. $R_2$) ensures that one can define a dg-functor $T_p:D_p\bt D_p\to D_{p+1}$  
(resp. $\inh _p:D_p^{op} \bt D_p\to D_{p+1}$) by 
$f_1\otimes f_2\mapsto T'_p(f_1,f_2)$ (resp. $f_1^{op}\otimes f_2\mapsto \inh' _p(f_1,f_2)$).
$R_3$ ensures that $a_{x_1,x_2,x_3}$, $\tau_{x_1,x_2}$ and $u_x$ are isomorphisms. \\

When $p\geq 0$, We put $\smor =\tens\sqcup \Int$ where $\tens$ and $\Int$ are defined by the same 
formula as above. To define $\srela$ we define four sets of relations $R_1$, $R_2$, $R'_3$, and $R'_4$. 
$R_1$ and $R_2$ are the ones defined above and 
\begin{eqnarray*}
R'_3&=&\{\:
T'_p(I_{p-1}(f_1), I_{p-1}(f_2))-T_{p-1}(f_1\otimes  f_2)\: \},\\
R'_4&=&\{\: 
\inh '_p(I_{p-1}(f_1), I_{p-1}(f_2))-\inh _{p-1}(f_1^{op}\otimes  f_2)\: \}.
\end{eqnarray*}
We set $\srela =R_1\sqcup R_2\sqcup R'_3\sqcup R'_4$ and $D_{p+1}:=D_p[\sob ,\smor ]/\srela$.
$R'_3$ and $R'_4$ ensure the  compatibility involving $I_{-}$, $T_{-}$ and $\inh _{-}$ 
and we have completed the induction. \qedsymbol \\

Now, we shall define the free closed tensor dgc $\fcl (C)$ associated to $C$. Set
\begin{eqnarray*}
D'&=&\colim _j\, (I_j,D_j) ,\\
T'&=&\colim _j\, T_j:D'\bt D'\to D', \\
\inh '&=&\colim _j\, \inh _j :D'^{op}\bt D'\to D',
\end{eqnarray*} 
where we regard $D'\bt D'\cong \colim _j(D_j\bt D_j)$ and $D'^{op}\bt D'\cong \colim _j
(D_j^{op}\bt D_j)$. These are well-defined by the compatibility of $I_j$'s with $T_j$'s 
and of $I_j$'s with $\inh _j$'s. Let $J$ be the $T' $-closed 
and 
$\inh'$-closed ideal generated by the relations which ensure the following conditions.
\begin{itemize}
\item $a$, $\tau$ and $u$ are natural isomorphisms and all of  
the coherence diagrams required in 
the definition of closed tensor dgc are commutative.

\item The morphisms of complex 
$\phi _{x,y,z}:\homo_{D'}(x\otimes y,z)\rightarrow \homo_{D'}(x,\inhom(y,z))$
given by $f\mapsto \inh' (id_y\otimes f)\circ coev_x^y $ form a natural isomorphism 
whose inverse $\varphi_{x,y,z}:\homo_{D'}(x,\inhom (y,z))\rightarrow 
\homo_{D'}(x\otimes y,z)$ is given by $g\mapsto \ev_y^z\circ T'(g\otimes id_y)$. 

\item $(x_0\oplus x_1,s_i^{x_0,x_1},p_i^{x_0,x_1})$ is a biproduct of $x_0$, $x_1$ (see \cite[P.190]{mac}).
\end{itemize}

Obviously these relations are represented by elements of $\mor (D')$. We set 
\[
\fcl (C):=D'/J.
\]
Using the universality of each $D_j$, one can check the functor $:\dgc\to\dgcl$ 
given by $C\mapsto \fcl (C)$ is a left adjoint 
of the forgetful functor and we have completed the construction of the functor $\fcl$ \qedsymbol

\subsection{A model category structure on $\dgcl$}\label{proofofmodel}
\subsubsection{Limit and colimit}
 We must show the category $\dgcl$ is closed under small limits and colimits. 
Limits are equal to those of underlying dg-categories with the additional structures naturally defined 
on them. As for colimits, pushouts can be constructed by induction. Let $C_1\leftarrow C_0
\rightarrow C_2$ be a diagram in $\dgcl$. We put $D_{-1}:=C_1\sqcup _{C_0}C_2$, the pushout in $\dgc$ (see 
\cite{tab1})  
and attach objects and morphisms step by step similarly  to the construction of $\fcl$ using universal dgc's. 
Infinite coproducts are similar.
\subsubsection{Functorial path object}\label{pathobj}
For the path object argument (Theorem \ref{pathobjarg}), we need a functorial path object in $\dgcl$, that is, a pair of 
\begin{itemize}
\item an endo-functor $\mca{P}:\dgcl\to \dgcl$ and  
\item a sequence of natural transformations 
$\{\,C\to \mca{P}(C)\stackrel{d_0\times d_1}{\to }C\times C \,\}_{C\in \dgcl}$ 
which is a factorization of the diagonal, such that $\mca{U}(C)\to \mca{U}(\mca{P}(C))\to 
\mca{U}(C\times C)$, where $\mca{U}:\dgcl\to\dgcu $ is the forgetfull functor, 
is a path object diagram in $\dgcu$ for any $C\in\dgcl$.\\
\end{itemize}

\indent Note that $k\to \nabla (1,*) \stackrel{d_0\times d_1}{\to }k\times k$ is a path object 
in the category of commutative dg-algebras over $k$. (For the notations, see subsection \ref{scdga}.) 
We define $\mca{P}(C)$ as follows.
\begin{itemize}
\item An object of $\mca{P}(C)$ is an isomorphism in $Z^0C$.
\item For isomorphisms $f:c_0\to c_1$, $f':c_0'\to c_1'\in Z^0C$ 
\[
\homo _{\mca{P}(C)}(f,f'):=\homo _C(c_0,c_0')\ot_k\nabla (1,*)
\]
and the composition is given by 
$(\beta \otimes\eta )\circ (\alpha \otimes \omega )
:=(-1)^{\deg \eta\cdot \deg \alpha}(\beta \circ \alpha)\otimes (\eta\cdot \omega)$, 
where $\alpha \in \homo _C(c_0,c_0')$, $\beta\in \homo _C(c_0',c_0'')$ and $\omega ,\eta \in \nabla (1,*)$. 
\item The additional structures are defined by those of $C$. For example, 
\begin{eqnarray*}
\inhom _{\mca{P}(C)}(f,f')&:=&\inhom _C(f^{-1},f'):\inhom (c_0,c_0')\to \inhom (c_1,c_1')\\
\inhom _{\mca{P}(C)}(\alpha \otimes \omega ,\beta\otimes \eta )&
:=&(-1)^{\deg \omega\cdot \deg \beta}\inhom (\alpha ,\beta )\otimes (\omega \cdot \eta ), 
\end{eqnarray*}
where $\alpha \in \homo _C(c_0,c_0')$, $\beta\in \homo _C(d_0,d_0')$ and $\omega ,\eta \in \nabla (1,*)$. 
\end{itemize}
$d_0:\mca{P}(C)\to C$ is given by $(c_0\to c_1)\mapsto c_0$ and $id\ot _kd_0:\homo (c_0,c_0')\ot\nabla (1,*)
\to \homo (c_0,c_0')$, $d_1:\mca{P}(C)\to C$ is given by $(c_0\to c_1)\mapsto c_1$ 
and $(f'_*\circ (f^{-1})^*)\otimes d_1:\homo (c_0,c_0')\ot\nabla (1,*)
\to \homo (c_1,c_1')$ and $C\to \mca{P}(C)$ by 
$c\mapsto (c\stackrel{id}{\to }c)$.
\begin{lem}
 With above notations, the sequence
\[
\xymatrix{\mca{U}C\ar[r]^{\mca{U}i}&\mca{UP}C\ar[rr]^{\mca{U}d_0\times\mca{U}d_1\ \ }&&\mca{U}C\times\mca{U}C}
\]
is a path object of $\mca{U}C$ in $\dgcu$.
\end{lem}
\begin{proof}
We shall show $\mca{U}i$ is a quasi-equivalence and $\mca{U}d_0\times \mca{U}d_1$ is a fibration in $\dgcu$. In the following we omit $\mca{U}$.
Clearly any object $f:c_0\to c_1\in \mca{P}C$ is isomorphic to $id_{c_0}$ and $k\to \nabla (1,*) \stackrel{d_0\times d_1}{\to }k\times k$ is a path object in the category of commutative dg-algebras with the usual model structure so we see that $i$ is a quasi-equivalence and $d_0\times d_1$ induces surjections on complexes of morphisms. Let  $f:c_0\to c_1\in \mca{P}C$ be an object and $(g,g'):(c_0,c_1)\to(c,c')\in C\times C$ be an isomorphism in $H^0(C\times C)=Z^0(C\times C)$. Note that $(c,c')=(d_0\times d_1)(g'\circ f\circ g^{-1})$, where $g'\circ f\circ g^{-1}$ is considered as an object of $\mca{P}C$.  Consider $g$ as an isomorphism $f\to g'\circ f\circ g^{-1}\in Z^0\mca{P}C$ via the injection $\homo_C(c_0,c)\to \homo_C(c_0,c)\otimes_k\nabla(1,*)=\homo_{\mca{P}C}(f,g'\circ f\circ g^{-1})$. 
It is clear that $(d_0\times d_1)g=(g,g')$ so by above assertion, $d_0\times d_1$ is a fibration in $\dgcu$. \end{proof}
Now one can prove the following theorem 
using   the free functor 
$\fcl ^{\geq  0}:\dgcu\to \dgcl$ and the functorial path object $\mca{P}$.
\begin{thm}\label{modelstrondgcl}
The category $\dgcl$ has a 
cofibrantly generated model category structure where 
weak equivalences and fibrations are
defined as follows.
\begin{itemize}
\item A morphism $F:C\rightarrow D\in\dgcl$ is a weak equivalence if and only if
it is a quasi-equivalence.
\item A morphism $F:C\rightarrow D\in\dgcl$ is a fibration if and only if it satisfies 
the following two conditions.
\begin{itemize}
\item For $c,c'\in \ob (C)$ the morphism $F_{(c,c')}:\homo _C(c,c')\to \homo _D(Fc,Fc')$ is an epimorphism. 

\item For any $c\in \ob (C)$ and any isomorphism $f:Fc\to d'\in Z^0(D)$, there exists an isomorphism 
$g:c\to c' \in Z^0(C)$ such that $Z^0(F)(g)=f$.
\end{itemize}
\end{itemize}

\end{thm}

\begin{proof} Apply Theorem \ref{pathobjarg} to the forgetful functor $\mca{U}:\dgcl\longrightarrow\dgcu$. Note that the fibrations of the statement correspond to those of Theorem \ref{dgcu} as $H^0C=Z^0C$ for $C\in\dgcl$.\end{proof}

\indent We also consider the augmented category. Let $\vect '$ denote the category of all finite dimensional $k$-vector spaces and $k$-linear maps.
We let $\mb{1}=k$ regarded as a $k$-vector space and fix a 0-dimensional vector space $\mb{0}$. With these 
distinguished objects,  
$\vect '$ has a closed tensor structure with the 
the usual operations $\otimes$, $\inhom$ and $\oplus$. 
We denote by $\vect$ the full subcategory of $\vect'$ consisting of objects 
represented by  words generated by 
${\mb{1},\mb{0}}$ with the operations $\otimes$, $\inhom$ and $\oplus$. $\vect$ is small and 
 it is (isomorphic to) 
an initial object of $\dgcl$. 
\begin{defi}\label{defofdgclp}
We call the over category $\dgcl /\vect$ the \textup{category of closed tensor dg-categories with fiber functors} and denote it by $\dgclp$. For an object 
$C=(C,\omega_C)\in\dgclp$ we call $\omega_C :C\to\vect$ the \textup{fiber functor of $C$}. 
\end{defi}
The following is a corollary of Theorem \ref{modelstrondgcl} (see \ref{modelcategories}).
\begin{cor}\label{modelstrondgclp}
$\dgclp$ has a model category structure induced by that of $\dgcl$.
\end{cor}

\section{The Sullivan-de Rham  equivalence for finite fundamental group}
\indent The purpose of this section is  to prove an equivalence theorem for the homotopy category of spaces whose fundamental group is finite and whose homotopy groups are finite dimensional $\bb{Q}$-vector spaces. We call the equivalence the Sullivan-de Rham equivalence because it is a direct generalization of the one in \cite{gug}. \\

For an abstract group $G$,  $\repg$ stands for 
the category of finite dimensional $k$-linear representations of $G$ whose underlying vector spaces belong to $\vect$.
$\repg$ has a   closed tensor  structure such that the forgetful functor $\omega _G:\repg\to \vect$ is a morphism of 
closed tensor  categories. 
\subsection{The generalized de Rham functor}\label{scdga}
\indent In this subsection, we define a Quillen pair 
\[
\xymatrix{
\tdr:\ssetp \ar@<3pt>[r]& (\dgclp)^{op}:\mathit{Sp}.\ar@<3pt>[l]}
\] 
We first recall the notion of standard simplicial commutative
dga $\nabla (*,*)$ over $k$ from \cite[Section 1]{gug}. 
Let $p\geq 0$ and $\nabla  (p,*)$ be the commutative graded algebra over $k$ 
generated by indeterminates $t_0,\dots,t_p$ of degree 0 and $dt_0,\dots,dt_p$
of degree 1 with relations
\[
t_0+\cdots +t_p =1, \ \ \
dt_0+\cdots +dt_p=0.
\]

We regard $\atp$ as a cdga with the differential given by $d(t_i):=dt_i$. We
can define simplicial operators
\[
d_i:\atp \rightarrow \nabla (p-1,*),\ \ \  
s_i:\atp \rightarrow \nabla (p+1,*),\ \ 0\leq i\leq p 
\]
(see \cite{gug}) and we also regard $\at$ as a simplicial commutative dga.\\
\indent The following definition is adopted in \cite{nsrat}
\begin{defi}
Let $\vect ^{iso}$ be the subcategory of $\vect$ consisting of
all objects and isomorphisms. Let $K$ be a simplicial set.\\
\indent A \textup{local system} $\loc $ \textup{on} $K$ is a functor 
$(\Delta K)^{op}\rightarrow \vect ^{iso}$ such that
for any simplex $\sigma\in \Delta K$, any degeneracy operator $s_i$ and
 any morphism $f:s_i\sigma\rightarrow \sigma$, 
$\msc{L} (\sigma)=\msc{L} (s_i(\sigma ))$ and $\msc{L} (f)=id_{\loc (\sigma)}$.\\
\indent A \textup{morphism of local systems} $\msc{L}\rightarrow \msc{L}'$ is a 
natural transformation $I \circ \msc{L} \Rightarrow 
I\circ \loc':(\Delta K)^{op}\rightarrow  \vect$, where 
$I:\vect ^{iso}\rightarrow \vect$ is natural inclusion functor. We define the
 \textup{tensor} $\loc\otimes \loc'$, the \textup{internal 
hom object} $\mf{Hom}(\loc,\loc')$ and the \textup{coproduct} $\loc \oplus \loc'$
 of two local systems $\loc$, $\loc'$ objectwisely by using those of $\vect $.
We denote by $\lock$ be the closed tensor category of local systems on $K$. 
If $K$ is pointed, $\lock$ is regarded as a closed tensor  category with the fiber functor 
given by the evaluation at the base point.
\end{defi}
It is well-known that for pointed connected simplicial set $K$, there exists an equivalence of closed tensor  categories $\lock\stackrel{\sim}{\to}\mathit{Rep}(\pi _1(K))$ which is functorial in $K$. In the following, we sometimes identify local systems with representations of the fundamental group, fixing such an equivalence.
\begin{defi}
Let $K$  be a simplicial set and $\loc$ be a local system on $K$. The \textup{de Rham complex of $\loc$-valued 
polynomial forms $A_{dR}(K,\loc)\in \ms{C}(k)$} is defined as follows. For each $q\geq 0$, the degree $q$ part is given by  
\[
A_{dR}^q(K,\loc )=\lim\nolimits_{\Delta K^{op}}\atq\ot _k\loc .
\]
Here $\atq$ is regarded as a functor from $\Delta K^{op}$ to the category of $k$-vector spaces by composed with the functor 
$\Delta K^{op}\rightarrow \Delta^{op}$, the limit is taken in the 
category of possibly infinite dimensional $k$-vector spaces. For $q\leq -1$, we set $A_{dR}^q(K,\loc )=0$. The differential $d:A_{dR}^q(K,\loc )\to A_{dR}^{q+1}(K,\loc )$ is defined from
 the differential $d:\nabla(*,q)\to \nabla(*,q+1)$.
\end{defi}
 We shall define the generalized de Rham functor
\[
\tdr:\sset\longrightarrow (\dgcl )^{op}.
\]
This  
is a natural generalization of the de Rham functor of \cite[Definition 2.1]{gug}. 
For $K\in\sset$ we define 
a closed tensor dgc $T_{dR}(K)$ as follows. An object is a local system on $K$ and 
$Hom_{\tdrk}(\loc,\loc ')=A_{dR}(K,\mf{Hom}(\loc,\loc '))$. The composition is defined 
from that of the category of vector spaces and the multiplication of $\at$, i.e.,
\[
(\eta \cdot b)\circ (\omega \cdot a):=(\eta\cdot\omega )\cdot (b\circ a)
\]
for $\omega,\eta\in\at$, $a\in \mf{Hom}(\loc,\loc')$ and $b\in \mf{Hom}(\loc ',\loc '')$. 
The additional structures $\otimes$, $\inhom$ and $\oplus$ are defined similarly. (We agree that)  
$\tdr (\emptyset )$ is a terminal object of $\dgcl$ and we always identify $\tdr (*)$ with 
$\vect $ via the isomorphism $\loc \mapsto \loc (*)$. 
For each morphism $f:K\rightarrow L\in\sset$ we associate a morphism of closed tensor dgc's 
$f^*:\tdrl\rightarrow \tdrk$ by 
\[
(\Delta (L)^{op}\stackrel{\loc}{\to}\vect ^{iso})\longmapsto 
(\Delta (K)^{op}\stackrel{\Delta f^{op}}{\to}\Delta (L)^{op}\stackrel{\loc}{\to}\vect ^{iso}).
\]
Thus we have defined a functor 
$\tdr:\sset\longrightarrow (\dgcl )^{op}$. \\ 
\indent Let $C\in\dgcl$. We define a functor $\mathit{Sp}:(\dgcl)^{op}\rightarrow \sset$ 
by
\[
\mathit{Sp}(C)_n=\homo _{\dgcl}(C,\tdr (\Delta ^n))
\]
with obvious simplicial operators.
\begin{lem}\label{quillenadjoint}
The functor $\mathit{Sp}$ is a right adjoint of $\tdr$ and the adjoint pair 
is a Quillen pair between $\sset$ and $(\dgcl )^{op}$. So it induces a Quillen pair between pointed categories:

\[
\xymatrix{
\tdr:\ssetp \ar@<3pt>[r]& (\dgclp)^{op}:\mathit{Sp}.\ar@<3pt>[l]}
\]\end{lem}
\begin{proof} The first assertion is clear. To show the second one, it is 
enough to examine the generating cofibrations and trivial cofibrations of $\sset$. 
One can check easily the condition about lifting of isomorphisms and the proof
 reduces to the case of constant coefficients, see \cite[Section 1, 2]{gug}. 
The third follows from the second.
\end{proof}

\subsection{Tannakian dg-categories and equivariant commutative dg-algebras}\label{tannandeqalg}
\indent In this subsection, we introduce the category of Tannakian dg-categories of finite type, which we will prove corresponds to the category $\ssetpfq$ via $\tdr$, and 
compare it with the category of equivariant dg-algebras. Throughout 
this subsection, we assume $k=\bb{Q}$.\\
\subsubsection{Tannakian dg-categories of finite type}
\begin{defi}\label{tannakian}
Let $(C,\omega_C )\in \dgclp$. We say $(C,\omega _C)$ is a \textup{Tannakian dg-category of finite type} if 
the following conditions are satisfied.
\begin{itemize} 
\item $(Z^0(C),Z^0(\omega ))$ is equivalent to $(\repg ,\omega_G)$, the closed tensor  category of finite dimensional representations of $G$ with the forgetful functor.  More precisely, there exists a finite chain of morphisms of closed tensor  categories with fiber functors which are equivalences of underlying categories:
\[
(Z^0(C),Z^0(\omega ))\stackrel{\sim}{\to}(T_1,\omega _1)\stackrel{\sim}{\leftarrow}\cdots \stackrel{\sim}{\to}(T_n,\omega _n)\stackrel{\sim}{\leftarrow}(\repg ,\omega_G).
\]
\item For each $c_0,c_1\in \ob C$, $H^1(\homo _C(c_0,c_1))=0$ and $H^i(\homo _C(c_0,c_1))$ is 
finite dimensional for  $i\geq 2$

\end{itemize}
We denote by $\tannf$ 
the full subcategory of $\dgclp $ consisting of Tannakian dgc's of finite type.
\end{defi}
Clearly $\tannf$ is stable under weak equivalences of $\dgclp$.\\
\indent If we use Tannakian theory, we get an  internal characterization of the subcategory $\tannf$ (which is not used in the rest of the paper, see \cite[Theorem 2.11,Proposition 2.20 (a)]{dmos}):  
\begin{prop}
 An object $(C,\omega_C )\in \dgclp$ belongs to $\tannf$ if and only if the following conditions are satisfied.\begin{itemize} 
\item The additive category $Z^0(C)$ is an abelian category and the functor $Z^0\omega_{C}:Z^0(C)\to \vect$ is exact and faithful. 
\item $\homo_{Z^0C}(\mb{1},\mb{1})\cong k$.
\item There exists an object $c\in Z^0C$ such that any object of $Z^0(C)$ is a sub-object of a finite coproduct $c^{\oplus N}$ for some $N$.
\item For each $c_0,c_1\in \ob C$, $H^1(\homo _C(c_0,c_1))=0$ and $H^i(\homo _C(c_0,c_1))$ is 
finite dimensional for  $i\geq 2$.
\end{itemize}
\end{prop}

\subsubsection{Equivariant commutative dg-algebras}
\indent Let $G$ be an abstract group. Let $\modg$ be the category of 
possibly infinite dimensional right $G$-modules over 
$k$ and $\dgmodg$ be the category of cochain complexes over $\modg$. $\dgmodg$ has a structure of 
symmetric monoidal category as usual and we denote by $\dgalgg$ the category of commutative monoids over 
$\dgmodg$. We call an object of $\dgalgg$ a $G$-equivariant commutative dg-algebra, in short, 
$G$-cdga. \\
\indent $\dgmodg$ has a model category structure such that a morphism is a weak equivalence 
(resp. fibration) if and only if it is a quasi-isomorphism (resp. levelwise epimorphism). The 
following lemma follows from the path object argument.
\begin{lem}
$\dgalgg$ has a model structure such that a morphism is a weak equivalence 
(resp. a fibration) if and only if it is a quasi-isomorphism (resp. levelwise epimorphism).
\end{lem}

\begin{defi}\label{defofeqalg}
\textup{(1)} The \textup{category of equivariant cdga's} $\eqalg$ is defined as follows.
\begin{itemize}
\item An \textup{object} is a pair $(G,A)$ of a group $G$ and a $G$-cdga $A$. 
\item A \textup{morphism} $f:(G,A)\rightarrow (H,B)$ is a pair of a group homomorphism 
$f^{gr}:H\rightarrow G$ and a morphism of $H$-cdga $f:(f^{gr})^*A
\rightarrow B$.
\end{itemize} 
We say an equivariant cdga $(G,A)$ is \textup{of finite type} if $G$ is a finite group 
and $H^iA$ is finite dimensional for any $i$ and \textup{1-connected} if $H^0A\cong k$ 
and $H^1A\cong 0$. we denote by $\eqalgf$ the full subcategory of $\eqalg$ consisting 
of 1-connected cdga's of finite type and by $\eqalgfp$ the over category $\eqalgf /(e,k)$, 
where $e$ is a trivial group.\\
\textup{(2)} We define the \textup{homotopy category} $\ho (\eqalgfp)$ as the localization of $\eqalgfp$ obtained by inverting all the maps whose group homomorphisms are isomorphisms and whose cdga homomorphisms are quasi-isomorphisms. \\
\textup{(3)} Let $f_1,f_2:(G,A)\to (H,B)$ be two morphisms of $\eqalgfp$. we say $f_1$ and $f_2$ are \textup{right  homotopic}, written $\sim_r$  if $f_1^{gr}=f_2^{gr}$ and if  $f_1$, $f_2:(f_1^{gr})^*A\to B$ are right  homotopic as morphisms of $\dgalgh /k$ with respect to its model structure (see \cite[DEFINITION 1.2.4]{hov}).
\end{defi}
The following lemma is proved by arguments similar to the proofs of \cite[Proposition 1.2.5, Theorem 1.2.10]{hov}. 
\begin{lem}\label{htpyineqalg}
Let $(G,A), (H,B)\in\eqalgfp$ and suppose $A$ is cofibrant as an object of $\dgalgg$. Then  the  relation $\sim _r$ on $\homo_{\eqalgfp}((G,A),(H,B))$ is an equivalence relation and there exists a bijection 
\[
\homo _{\eqalgfp}((G,A),(H,B))/\sim_r\cong \homo _{\ho (\eqalgfp)}((G,A),(H,B))
\]
which takes a class represented by a map $f:(G,A)\to (H,B)$ to the image of $f$ by the canonical functor $\eqalgfp\to \ho (\eqalgfp )$.
\end{lem}

\subsubsection{Comparison}
We shall define two functors 
\[
\ms{T}, \ms{T}^c:\eqalgfp \longrightarrow \dgclp.
\]
Let $A=(G,A)\in \eqalgfp$. 
For $V\in\ob(\mathit{Rep}(G))$ we 
define a complex $A\ot ^GV$ by
\[
A\ot ^GV:=\{\: \Sigma _ja_j\ot v_j\in A\ot _kV\:|\:\Sigma _j(a_j\cdot g)\ot v_j=\Sigma _ja_j\ot (g\cdot v_j) 
\textrm{ for $\forall g\in G$.}\:\}
\]
and set 
\[
\ob(\ms{T}A):=\ob(\mathit{Rep}(G)),\ \ \ \ \homo _{\ms{T}A}(V,W):=A\ot^G\inhom_{\repg}(V,W),
\]
where $\inhom_{\repg}$ is the internal hom of $\repg$.  We define composition, closed tensor structure of $\ms{T}(A)$ using corresponding structures of $\repg$ and the multiplication of $A$. A morphism $f:(G,A)\rightarrow (H,B)$ of $\eqalgfp$ gives a functor 
$(f^{gr})^*:\repg \rightarrow \reph $ so $f$ induces a morphism 
$\ms{T}f:\ms{T}A\rightarrow \ms{T}B$ of closed tensor dg-categories.  The augmentation $A\to k$ defines a fiber functor $\ms{T}A\rightarrow \ms{T}k\cong \vect$ and  we have defined a functor $\ms{T}:\eqalgfp \rightarrow \dgclp$.

\begin{exa}
Let $(G,k)\in\eqalgfp$ denote the equivariant cdga whose underlying cdga is $k$ and whose group action is trivial. $\ms{T}(G,k)$ is isomorphic to $(\repg, \omega _G)$.
\end{exa}

\indent $\ms T^c$ is defined as follows. Let $\wcl(\ob(\mathit{Rep}(G)))$ be the set of words freely generated by 
$\ob(\mathit{Rep}(G))\sqcup\{\mb{1},\mb{0}\}$ 
with operations $\widehat{\otimes}$, $\widehat{\oplus}$ and $\widehat{\inhom} $ (see sub-subsection \ref{freecl}). Let $R:\wcl( \ob(\mathit{Rep}(G)))\rightarrow\ob(\mathit{Rep}(G))$ 
be the  function defined inductively, by 
\begin{itemize}
\item $RX=X$ for $X\in \ob(\mathit{Rep}(G))$, 
\item $R(X\widehat{\otimes} Y)=(RX)\otimes (RY)$,  
$R(\widehat{\inhom} (X,Y))=\inhom (RX,RY)$, and $R(X\widehat{\oplus} Y)=(RX)\oplus (RY)$
\end{itemize}
We define $\ms{T}^c$ as the "pullback" of $\ms{T}$ by $R$. Precisely, we set 
\[
\ob\ms{T}^cA:=\wcl(\ob(\mathit{Rep}(G))),\ \ \ \homo _{\ms{T}^cA}(X, Y):=\homo _{\ms{T}A}(RX,RY)
\] Also the closed tensor structure on $\ms{T}^cA$ is  defined by "pullback" by $R$.  
For example, the tensor structure is given by the operation $\widehat{\otimes}$ of $\wcl(\ob(\mathit{Rep}(G)))$. $R$ defines a morphism $R_A:\ms{T}^cA\to \ms{T}A$ and we define the fiber functor of $\ms{T}^cA$ by the composition 
$\ms{T}^cA\stackrel{R_A}{\to}\ms{T}A\stackrel{\omega _{\ms{T}A}}{\to}\vect$. 
Thus, we have defined a functor $\ms{T}^c:\eqalgfp \rightarrow \dgclp$.

Obviously,  $\ms{T}^cA$ is naturally equivalent to $\ms{T}A$ via $R_A$ and $\ms{T}A$ is simpler but  $\ms{T}^cA$ is convenient to model categorical arguments  because  of the following:
\begin{lem}\label{cofibrancy}
For a finite group $G$, $\T^c(G,k)$ is cofibrant in $\dgclp$. 
\end{lem}
\begin{proof}
 For closed tensor  categories $C,D$, consider the lifting problem
\[
\xymatrix{
&C\ar[d]^P\\
\T^c(G,k)\ar@{-->}[ur]\ar[r]^F&D,}
\]
where $F$ and $P$ are morphisms of closed tensor  categories and $P$ is an equivalence of categories which induces a surjective map on the sets of objects. We can find a lifting $\widetilde F:\T^c(G,k)\to C$ as a functor but $\widetilde F$ may not be a morphism of closed tensor categories. For example, $\tF(X)\otimes \tF(Y)$ and $\tF(X\otimes Y)$ are isomorphic, but not always equal. We shall modify $\tF$ so that it becomes a morphism of closed tensor  categories. As  $\ob(\T^c(G,k))$ is freely generated by $\ob(\repg)\sqcup \{\mb{1},\mb{0}\}$, we can define a function $\tF':\ob(\T^c(G,k))\to \ob(C)$ such that $\tF'(X)=\tF(X)$ for $X\in \ob(\repg)$ and $\tF'$ preserves $\otimes, \inhom,\oplus, \mb{1}$ and $\mb{0}$. For a morphism $f\in\homo _{\T^c(G,k)}(X,Y)$, we define $\tF'(f)\in\homo _C(\tF'X,\tF'Y)$ as the composition
\[
\tF'X\stackrel{\varphi_X}{\to}\tF X\stackrel{\tF(f)}{\to}\tF Y\stackrel{\varphi_Y^{-1}}{\to}\tF'Y.
\]  
Here, $\varphi_Z:\tF'(Z)\to \tF(Z)$ is the unique isomorphism such that $P(\varphi _Z)=id_{F(Z)}$. Thus we have defined a functor $\tF':\T^c(G,k)\to C$. This is clearly a morphism of closed tensor categories and we have completed the proof.
\end{proof}
To see some basic properties of $\ms{T}$, we need a few elementary preparations. For a finite group $G$ let $\vgr $ be the vector spaces of $k$-valued functions on $G$. (This vector space does not belong to $\vect$ but it is isomorphic to an object of $\vect$ and so we fix an isomorphism and we deal with $\vgr$ as if it belonged to $\vect$ via the isomorphism.) There  are two actions $\rho$, $\varrho$ of $G$ on  $\vgr$. 
\[
[\rho (g)\alpha](g')=\alpha(g'g),\ \  [\varrho (g)\alpha](g')=\alpha(gg'),\ \ \  \alpha:G\to k\in \vgr.
\]
$\rho$ is a left action and $\varrho$ is a right action. We call the representation $(\vgr,\rho )$ the  regular representation of $G$ and sometimes we omit $\rho$. Let $V\in\ob(\repg )$. In the following we use the homomorphism 
\[
\phi _V:V\to \inhom _{\repg}(V_u^{\vee },\vgr ),\ \ \phi  _V(v)(v')(g)=\langle v',gv\rangle
\]
where $V_u$ is the trivial representation with the same underlying vector space as $V$. This is a monomorphism and has a retraction.
\begin{lem}\label{presqequi}
\textup{(1)} For each $(G,A)\in \eqalgfp$,  the morphism of complexes 
\[
A\longrightarrow A\otimes ^G\vgr, \ \ \   a\longmapsto \Sigma_{g\in G}a\cdot g\otimes \delta _g,
\] where $\delta _g$ is the $\delta$-function at $g\in G$, is an isomorphism.\\
\textup{(2)} Let $f:(G,A)\to (G,B)$ be a morphism of $\eqalgfp$ such that $f^{gr}=id$. 
If $f$ is a quasi-isomorphism of $G$-cdga's, 
$\ms{T}f:\ms{T}A\to \ms{T}B$ is a quasi-equivalence.\\
\textup{(3)} $\ms{T}A$ is a Tannakian dgc of finite type for any $A\in \eqalgfp$.  
\end{lem}
\begin{proof} (1) is easy. (2) is a consequence of (1) and the fact that for any $V\in \repg$, $A\otimes^GV$ is a retract of $A\otimes ^G\inhom _{\repg}(V_u^{\vee },\vgr )\cong A^{\oplus \dim V}$. (3) follows from the fact that 
$(G,A)$ is quasi-isomorphic to  a $G$-cdga $(G,M)$ such that $M^0\cong k$, $M^1\cong 0$.  
\end{proof}

  Let $(C,\omega _C)$, $(D,\omega_D)\in\dgclp$ and $F$, $F':(C,\omega _C)\to (D,\omega_D)$ be 
two morphisms. We say $F$ and $F'$ are 2-isomorphic if there exists a natural isomorphism 
$i_-:F\Rightarrow F'$ such that $i_{c_0\otimes c_1}=i_{c_0}\otimes i_{c_1}$ 
and $\omega _D(i_c)=id_{\omega_C(c)}$. 
If $C$ and $D$ are Tannakian dgc's of finite type, $i$ is unique if it exists. \\
\indent (1) of the following is a rewrite of  \cite[Proposition 2.8]{dmos} for the case of finite group and (2) follows from (1).
\begin{thm}[\cite{dmos}]\label{tdforfinite}
Let $G$ be a finite group. Let $\omega=\omega_{\ms{T}^c(G,k)}$. Let $\aut(\omega)$ be the group of tensor preserving automorphisms of $\omega$ i.e.,
\[
\aut(\omega)=\{\alpha:\omega\Rightarrow \omega |\alpha _{X\otimes Y}=\alpha _X\otimes \alpha _Y\}.
\]
\textup{(1)} The homomorphism $\phi _G: G\to \aut(\omega )$  defined by $\phi _G(g)_X=r_{RX}(g):\omega(X)\to \omega(X)$, where $r_{RX}$ is the action of $G$ endowed with the representation $RX$, is an isomorphism of groups\\
\textup{(2)} Let $H$ be another finite group. Set $\omega'=\omega_{\ms{T}^c(H,k)}$. Let $F:\ms{T}^c(G,k)\to\ms{T}^c(H,k)$ be a morphism of closed tensor category such that $\omega'\circ F=\omega.$ $F$ is 2-isomorphic to  $\ms{T}^c(f)$ for some $f:(G,k)\to (H,k)\in \eqalgfp$.
\end{thm}
 
\indent The following theorem says the category of Tannakian dgc's of finite type and the category of 1-connected augmented equivariant dg-algebras of finite type are essentially the same. In the proofs of this theorem and Lemma \ref{pushout} we need  internal hom functors $\inhom$.
\begin{thm}\label{eqalg}
\textup{(1)} The functor $\ms{T}^c:\eqalgfp\to\tannf$ is fully faithful up to 2-isomorphisms. 
More precisely, for any $A$, $B\in\eqalgfp$ and 
$F:\ms{T}^cA\to \ms{T}^cB\in\tannf$ there exists a unique morphism $f:A\to B\in\eqalgfp$ such that $\ms{T}^cf$ is 
2-isomorphic to $F$. \\
\textup{(2)} Any object of $\tannf$ is equivalent to $\ms{T}A$ (and $\ms{T}^cA$) for some $A\in \eqalgfp$.  \\
\textup{(3)} Let $f:(G,A)\to (G,B)\in\eqalgfp$ be a morphism such that $f^{gr}$ is the identity. 
Suppose $f$ is a cofibration as a morphism of $\dgalgg$. 
Then $\ms{T}^cf:\ms{T}^cA\to \ms{T}^cB$ is a cofibration in $\dgclp$. In particular, if $(G,A)\in\eqalgfp$ is cofibrant as an object of $\dgalgg$, $\ms{T}^c(G,A)$ is cofibrant in $\dgclp$. \\
\textup{(4)} $\ms{T}$ and $\ms{T}^c$ induces an equivalence of categories $\ho (\eqalgfp)\simeq\ho(\tannf)$. 
\end{thm}
\begin{proof}
(1) Let $A=(G,A)$ and $B=(H,B)$. We show surjectivity up to 2-isomorphisms. By Theorem \ref{tdforfinite}, $F$ is 2-isomorphic to 
a morphism $F'$ such that $Z^0(F'):\T^c(G,k)\to \T^c(H,k)$ is $\T^c(f^{gr})$ for some group homomorphism $f^{gr}:H\to G$. So 
we may replace $F$ by such $F'$.\\
\indent   If $Z^0(F)$ is fixed as above,  $F$ is determined by 
\[
F_{(\mb{1},\vgr )}:\homo _{\ms{T}^cA}(\mb{1},\vgr )\rightarrow \homo _{\ms{T}^cB}(\mb{1},(f^{gr})^*\vgr )
\]
Let $A':=\homo _{\ms{T}A}(\mb{1},\vgr )$. The right action $\varrho$ on $\vgr$ defines a right action of $G$ on $A$  and pointwise multiplication 
$\vgr \ot \vgr \to \vgr$ defines a structure of cdga on $A'$. Thus, the tensor structure on $\ms{T}^c(A)$ defines a $G$-cdga structure on $A'$ and the isomorphism of Lemma \ref{presqequi} $A\to A'$ is an isomorphism of $G$-cdga's. Composing $F_{(\mb{1},\vgr )}$ with the 
homomorphism 
\[
\homo _{\ms{T}^cB}(\mb{1},(f^{gr})^*\vgr )\to \homo _{\ms{T}^cB}(\mb{1},\vhr )=:B'
\] 
induced 
by $(f^{gr})^*:(f^{gr})^*\vgr\to\vhr \in \reph$, we get 
a morphism of augmented equivariant cdga's $f':A'\to B'$. Composing with the isomorphisms $A\cong A'$, $B\cong B'$ defined in Lemma \ref{presqequi}, we get $f:A\to B$ such that $\ms{T}^cf=F$. The 
injectivity is clear from the above argument.\\
\indent (2) Let $T\in \tannf$. We may assume $Z^0T=\repg $ 
where $G$ is a finite group.
Let $V\in\repg$. There is a $G$-bimodule structure on $\inhom (V_u^{\vee },\vgr )$  
determined by the right action on $\vgr$ and the left action on $V_u$. 
The image of $\phi _V$ is $\{\alpha\in \inhom (V_u^{\vee },\vgr )|g\cdot \alpha 
=\alpha\cdot g \textrm{ for } \forall g  \}$. So if we put $A:=\homo _T(\mb{1},\vgr )$, the 
homomorphism $A\ot _kV_u\cong \homo _T(\mb{1},V_u\ot\vgr)\cong 
\homo _T(\mb{1},\inhom (V_u^{\vee },\vgr ))$ induces an isomorphism of complexes 
$A\ot^GV\cong \homo _T(\mb{1},V)$ via $\phi_V$. The composition $A\ot^G\inhom (V,W)\cong 
\homo _T(\mb{1}, \inhom (V,W))\cong \homo _T(V,W)$ defines a morphism $\ms{T}A\to T$ of Tannakian dgc's. 
Various naturalities ensure this is well-defined and this is clearly an equivalence. \\
\indent (3) Let $\ms{T}A\rightarrowtail T\stackrel{\sim}{\twoheadrightarrow} \ms{T}B$ 
be a factorization in $\dgclp$. Put $B':=\homo _T(\mb{1},\vgr )$. We regard $B'$ as a $G$-cdga 
which has an augmentation $B'\rightarrow k$. 
By an argument similar to the proof of (2), 
one can take morphisms $i:A\to B'$, $p:B'\to B\in \dgalgg /k$ 
and $F:\ms{T}B'\to T\in \tannf$ such that 
the composition $\ms{T}A\stackrel{\ms{T}i}{\to} \ms{T}B'\to T'$ is equal to 
the morphism of the factorization$\ms{T}A\to T'$ and the composition $\ms{T}B'\to T'\to \ms{T}B$ is equal 
to $\ms{T}B'\stackrel{\ms{T}p}{\to} \ms{T}B$. Then one can easily see  
$F$ has a retraction which preserves the morphisms from $\ms{T}A$. So  $\ms{T}i:\ms{T}A\to \ms{T}B'$ is 
a cofibration in $\dgclp$. 
One can see $\ms{T}f:\ms{T}A\rightarrow \ms{T}B$ is a retract of $\ms{T}i$ by 
using a lift of the following diagram.
\[
\xymatrix{
A\ar[r]^i\ar[d]^f&B'\ar[d]\\
B\ar@{-->}[ur]\ar[r]^{id}&B}
\] 
So it is a cofibration in $\dgclp$. The latter claim follows from Lemma \ref{cofibrancy}.\\
\indent (4) We  show the statement about $\T^c$. Then the one about $\T$ follows from it. By (2) it is enough to show the map 
\[
\T^c_{(A,B)}:\homo _{\ho (\eqalgfp )}(A,B)\longrightarrow \homo _{\ho (\dgclp)}(\T^cA,\T^cB)
\]
 is a bijection for each $A,B\in\eqalgfp$. We may assume $A$ is cofibrant as $G$-cdga and $\T^cB$ is fibrant. Then, by  (3), Lemma \ref{htpyineqalg} and \cite[Theorem 1.2.10,(ii)]{hov},the above sets of morphisms in homotopy categories are identified with the sets of  right homotopy classes of morphisms.  It is easy to see that for two morphisms $f_1$, $f_2:(G,A)\to (H,B)\in \eqalgfp$ $f_1$ and $f_2$ are right homotopic if and only if $\ms{T}f_1$ and $\ms{T}f_2$ are right homotopic in $\dgclp$  and 2-isomorphic morphisms in $\dgclp$ are right homotopic (see sub-subsection \ref{pathobj}). Then the claim follows from (1). 
\end{proof}

The following lemma is used in next subsection.
\begin{lem}\label{pushout}
Let $(G,A)\in \eqalgfp $. Let $\tilde T$ be an object of  $\tannf$ 
defined by $\ob (\tilde T)=\ob (\mathit{Vect})$ 
and $\homo _{\tilde T}(V,W)=A\ot _k\homo _{\vect }(V,W)$. 
Then there is a commutative diagram in $\tannf $
\[
\begin{CD}
\ms{T}^c(G,k)@>>>\vect \\
@VVV@VVV \\
\ms{T}^c(G,A)@>>>\tilde T \\
\end{CD}
\]
where the left vertical morphism is induced by the unit $k\to A$ and the bottom horizontal morphism is given by inclusion $A\otimes^G\inhom _{\repg}(RX,RY)\subset A\otimes \homo_{\vect} (\omega (X),\omega (Y))$ for each $X,Y\in\ob (\T^c(G,A))$. This diagram is a pushout diagram in $\dgclp$ and  a homotopy pushout diagram.
\end{lem} 
\begin{proof} 
We show the first assertion. The second one follows from it, Lemma \ref{cofibrancy}, and  
Theorem \ref{eqalg}, (3) (see also \cite[Lemma 5.2.6]{hov}). Let 
\[
\begin{CD}
\ms{T}^c(G,k)@>>>\vect \\
@VVV@VVV \\
\ms{T}^c(G,A)@>F>>C \\
\end{CD}
\]
be a commutative diagram in $\dgclp $. We define a homomorphism of complexes
\[
\tilde F_{(\mb{1},\omega (\vgr )\otimes V)}:\homo _{\tilde T}(\mb{1},\omega (\vgr )\ot V)
\to \homo _C(\mb{1},\omega (\vgr )\ot V)
\]
for $V\in \ob (\vect)$. Let $(\Sigma _{g\in G}a_g\otimes \delta _g)\ot v\in
\homo _{\tilde T}(\mb{1},\omega (\vgr )\ot V)\cong (A\otimes \omega (\vgr ))\otimes V$ where 
$a_g\in A$ and $v\in V$. Then $(\Sigma _{g'\in G}a_g\cdot g'\otimes \delta _{g'})\otimes v$ 
is regarded as an element of $\homo _{\ms{T}A}(\mb{1},\vgr \ot V)$. We set 
\[
\tilde F((\Sigma _{g}a_g\ot \delta _g)\ot v):=
\Sigma _g f_g\circ F((\Sigma _{g'}a_g\cdot g'\ot \delta _{g'})\ot v)
\]
where $f_g:\omega (\vgr )\to \omega (\vgr )\in \vect$ is 
\[
f_g(\delta _h)=\left\{
\begin{array}{ll}
\delta _g &\textrm{if }h=e \\
0         &\textrm{otherwise }
\end{array}\right.
\]
In general, one defines $\tilde F_{(V,W)}:\homo _{\tilde T}(V,W)\to \homo _C(V,W)$ 
using the embedding $\inhom (V,W)\to \inhom ((\inhom (V,W))^{\vee},\omega (\vgr ))\cong 
\inhom(V,W)\otimes \omega (\vgr) $. One can easily check $\tilde {F}_{(V,W)}$'s defines functor 
$\tilde F:\tilde T\to C\in\dgclp $ and $\tilde F$ makes appropriate diagrams 
commutative. 
\end{proof}

\indent  Let  $L\in\ssetpfq$. We take the universal covering $\widetilde L\to L$. The polynomial de Rham algebra $A_{dR}(\widetilde L)$ has a natural action of $\pi_1(L)$ induced by the action on $\widetilde L$. The construction $L\mapsto  (\pi_1(L), A_{dR}(\widetilde L))$ defines a functor $\widetilde{A_{dR}}:\ssetpfq\to (\eqalgfp )^{op}$. (Here, $\widetilde L$ is taken functorially in $L$.) \\
\indent We shall compare two functors $\tdr$ and $\widetilde{A_{dR}}$. 
\begin{prop}\label{naturaltrans}
Let $\ms{S}$ be either $\ssetpf$ or $\ssetpfq$ (see \ref{modelcategories}). Consider the following diagram.
\[
\xymatrix{
&(\eqalgfp)^{op}\ar[d]^{\ms{T}}\\
\ms{S} \ar[ur]^{\widetilde{A_{dR}}}\ar[r]^{T_{dR}} &(\tannf)^{op} \\
}
\]
There exists a natural transformation $\Phi :\tdr\Rightarrow \ms{T}\circ \widetilde{A_{dR}}$ 
such that for each $L\in\ssetpfq$, $\Phi _L:\tdr (L)\to\ms{T}\circ \widetilde{A_{dR}}(L)$ is an equivalence of underlying categories. 
\end{prop}
\begin{proof} The projection $p:\widetilde L\to L$ defines a morphism $p^*:\tdr (L)\to\tdr (\widetilde L)\simeq \vect \otimes A_{dR}(\widetilde L)$, where the closed tensor dg-category $ \vect \otimes A_{dR}(\widetilde L)$ is given by $\ob(\vect \otimes A_{dR}(\widetilde L))=\ob(\vect)$ and $\homo _{\vect \otimes A_{dR}(\widetilde L)}(V,W)=\homo_{\vect}(V,W)\otimes _kA_{dR} (\widetilde L)$. For two representations $V,W\in \mathit{Rep}(\pi_1(L))$, the morphism 
\[
p^*_{(V,W)}:\homo _{\tdr(L)}(V,W)\to \homo _{\vect}(V,W)\otimes A_{dR} (\widetilde L)
\] 
is a monomorphism and its image is precisely $\inhom (V,W)\otimes^{\pi_1(L)}\widetilde{A_{dR}}(L)$ so $p^*$ defines the required natural transformation. 
\end{proof}
In view of this proposition, we can produce some examples.
\begin{exa}
Let $G$ be a finite group and $L$ be a $K(G,1)$-space.  $\widetilde{A_{dR}}(L)$ is quasi-isomorphic to $(G,k)$ so $\tdr (L)$ is quasi-equivalent to $\ms{T}(G,k)\cong \repg$. 
\end{exa}

\begin{exa}
Let $L$ be the 2-dimensional real projective space $\bb{R}P^2$. $\widetilde{A_{dR}}(L)$ is quasi-isomorphic to $(\bb{Z}/2, M)$ where $M$ is a cdga freely generated by two generators $t$, $s$ with $\deg t=2$, $\deg s=3$ as   a commutative graded algebra and the differential is given by $d(t)=0$, $d(s)=t^2$. $\bb{Z}/2$ acts on $M$ by $g\cdot t=-t$ and $g\cdot s=s$ ($g$ is the generator). $\tdr (L)$ is quasi-equivalent to $\ms{T}(\bb{Z}/2,M)$. By definition, $Z^0(\ms{T}(\bb{Z}/2,M))\cong \mathit{Rep}(\bb{Z}/2)$. We shall describe the complex of morphisms $\homo_{\ms{T}(\bb{Z}/2, M)}(\mb{1},V)$ for each irreducible representation $V$. Let $V_-$ be the 1-dimensional representation where $g$ acts by the multiplication of $-1$.\\
\begin{tabular}{r|cccccccc}
                                      &0        &1      &2       &3         &4                       &5       &6 &7       \\
   \hline                                   
$\homo _{\ms{T}(\bb{Z}/2,M)}(\mb{1},\mb{1})=$&$\bb{Q}$&$0$       &$0$      &$\bb{Q} s$&$\bb{Q} t^2$&$0$& $0$  &$\bb{Q}st^2$               \\
$\homo _{\ms{T}(\bb{Z}/2,M)}(\mb{1},V_{-})=$   &$0$        &$0$      &$\bb{Q}t$&$0$        &$0$                      &$\bb{Q} st$  &$\bb{Q}t^3$ &$0$  \\
$M=$&$\bb{Q}$&$0$&$\bb{Q}t$&$\bb{Q} s$&$\bb{Q} t^2$&$\bb{Q} st$&$\bb{Q}t^3$&$\bb{Q}st^2$

\end{tabular}
\end{exa}

\subsection{The Sullivan-de Rham equivalence theorem for finite fundamental group}\label{suldr}

\indent In this subsection, we prove Theorem \ref{result2}. \\

\indent We first modify the right adjoint $\mathit{Sp}:(\dgclp)^{op}\longrightarrow \ssetp$. Let $\ssetpc$ be the full subcategory of $\ssetp$ consisting of connected pointed simplicial sets. 
The author cannot prove the image of $\ho (\tannf )$ by the functor 
$\bb{R}\mathit{Sp}:\ho (\dgclp)^{op}\longrightarrow\ho (\ssetp)$ is contained in $\ho (\ssetpc )$. We define a functor 
$\mathit{Sp}_0:(\dgclp)^{op}\to\ssetpc$ by saying that   
$\mathit{Sp}_0(C)$ is the connected component of $\mathit{Sp}(C)$ containing the base point 
for each $C\in\dgclp$. There is an obvious adjunction
\[
\xymatrix{
\tdr:\ssetpc \ar@<3pt>[r]& (\dgclp)^{op}:\mathit{Sp}_0.\ar@<3pt>[l]}
\]
This gives derived adjunction 
$\xymatrix{
\tdr:\ho(\ssetpc) \ar@<3pt>[r]&\ho (\dgclp)^{op}:\bb{R}\mathit{Sp}_0\ar@<3pt>[l]}$.

\begin{thm}\label{mainthm}
Suppose $k=\bb{Q}$. \\
\textup{(1)} The left Quillen functor $\tdr :\ssetp \to (\dgclp )^{op}$ induces an equivalence between homotopy categories:
\[
\xymatrix{\ho (\ssetpfq )\ar[r]^{\sim}&\ho (\tannf)^{op}}.
\]
\textup{(2)} Let $K\in\ssetpf$. The adjunction map 
\[
K\longrightarrow \bb{R}\mathit{Sp}_0\tdrk.
\]
is a fiberwise rationalization.
\end{thm}
To show this theorem, we need the following lemma and corollary.
\begin{lem}\label{kpi1}
Let $G$ be a finite group and $K$ be a $K(G,1)$-space. The unit of the adjunction 
$K\to \bb{R}\mathit{Sp}_0 \tdrk$ is a weak equivalence of simplicial sets.
\end{lem}
\begin{proof} $\tilde K$ is contractible so $\tdrk$ is quasi-equivalent to $\ms{T}^c(G,k)$. As $\ms{T}^c(G,k)$ is cofibrant, the 
morphism $K\to \bb{R}\mathit{Sp}_0 \tdrk$ is weak equivalent to $K\to \mathit{Sp}_0\ms{T}(G,k)$ which is  
the adjoint of the composition 
\[
\ms{T}^c(G,k)\stackrel{\sim}{\to}\repg\stackrel{\sim}{\to} \lock\cong Z^0\tdrk \to \tdrk .
\] 
One can see $\pi _i(\bb{R}\mathit{Sp}_0 \tdrk )=0$ for $i\geq 2$ by the adjunction. So it is enough to show 
$K\to \mathit{Sp}_0\ms{T}^c(G,k)$ gives an isomorphism of $\pi _1$. We may assume $K=N(G)$, 
the nerve of $G$. Both $K$ and 
$\mathit{Sp}_0\ms{T}^c(G,k)$ are fibrant, one can check this explicitly. A representative of a 
class in  
$\pi _1(\mathit{Sp}_0\ms{T}(^cG,k))$ is a morphism $F:\ms{T}^c(G,k)\to \tdr (\Delta ^1)$ of $\dgcl$ 
such that $d_0\circ F=d_1\circ F=\omega _{\ms{T}^c(G,k)}: \ms{T}^c(G,k)\to \vect$. We define an 
element $\alpha _F\in\mathit{Aut}^{\otimes}(\omega _{\ms{T}^c(G,k)})$ by $\alpha _F(V):=$
 the composition $F(V)_1\stackrel{(d_0)^{-1}}{\to}
F(V)_{\sigma}\stackrel{d_1}{\to}F(V)_0$ for $V\in\ob (\ms{T}^c(G,k))$ ($\sigma$ is the non-degenerate 1-simplex of $\Delta^1$).  
$\alpha _F$ corresponds to some $g\in G$ 
via canonical isomorphism $G\cong \mathit{Aut}^{\otimes}(\omega _{\ms{T}^c(G,k)})$ (see Theorem \ref{tdforfinite}). One can see 
$F$ represents the same class as $\mathit{Ev}_g:\ms{T}^c(G,k)\to \tdr (\Delta ^1)$, the evaluation 
at the edge corresponding to $g$ and the assertion follows.  
\end{proof}

\indent The following collorary follows from  Proposition \ref{naturaltrans} and  Lemma \ref{pushout}.
\begin{cor}\label{hfiber}
Let $L\in \ssetpfq$. Consider a homotopy fiber sequence 
\[
\widetilde L\longrightarrow L\longrightarrow K(\pi _1(L),1)
\]
where the right map induces isomorophism of $\pi _1$.
The corresponding sequence
\[
\tdr (K(\pi _1(L),1))\longrightarrow \tdr (L)\longrightarrow \tdr (\widetilde L)
\]
is a homotopy cofiber sequence in $\dgclp$. 
\end{cor}
\qedsymbol \\

\indent Now we shall prove Theorem \ref{mainthm}, (1). Let $L\in \ssetpfq$ and
$K$ be a 
 $K(\pi_1(L),1)$-space. Let $\widetilde L\to L\to K$ be a homotopy fiber sequence where the map $L\to K$ induces an isomorphism of $\pi _1$. 
Consider the following diagram.
\[
\begin{CD}
\widetilde L@>>>L@>>>K \\
@VVV@VVV@VVV \\
\rsp _0\tdr (\widetilde L)@>>>\rsp _0\tdr (L)@>>>\rsp _0\tdr (K) \\
\end{CD}
\] 
The left vertical arrow is a weak equivalence by original Sullivan's theory 
and the right vertical one is a weak equivalence by Lemma \ref{kpi1}. The 
bottom horizontal sequence is a homotopy fiber sequence by Corollary \ref{hfiber} and 
so is the top horizontal one by definition. 
Hence the middle vertical arrow is a weak equivalence. Thus $\tdr 
:\ho (\ssetpfq )\to \ho (\tannf )^{op}$ is fully faithful. Essential surjectivity 
follows from a similar argument and Theorem \ref{eqalg}. (2) follows from the above proof. \qedsymbol

\section*{Acknowledgements}
 I thank Masana Harada for informing me about 
non-Abelian Hodge theory, listening to my talks about this paper, and giving advices to 
improve readability of this paper. I also thank Yoshiyuki Kimura for letting me know 
utility of dg-category.

\end{document}